\documentclass[12pt]{article}

\usepackage{amsmath}
\usepackage{amssymb}
\usepackage{amsfonts}
\usepackage{eucal}
\usepackage{enumerate}
\usepackage{pgf}
\newcommand{\R}{\mathbb{R}}
\newcommand{\PP}{\mathbb{P}}

\newcommand{\C}{\mathbb{C}}
\newcommand{\K}{\mathcal{S}}

\newtheorem{theorem}{\bf Theorem}[section]
\newtheorem{proposition}[theorem]{\bf Proposition}

\newtheorem{lemma}[theorem]{\bf Lemma}
\newtheorem{definition}{\bf Definition}

\newtheorem{remark}[theorem]{\bf Remark}

\numberwithin{equation}{section}

\newenvironment{Proof}{\removelastskip\par\medskip
\noindent{\em Proof.}
\rm}{\penalty-20\null\hfill$\square$\par\medbreak}

\newenvironment{Proofy}{\removelastskip\par\medskip
\noindent{\em Proof} \hskip-0.08cm \rm}{\penalty-20\null\hfill$\square$\par\medbreak}

\def\div{{\mathrm{{\rm div \hskip0.048cm }}}}

\allowdisplaybreaks

\usepackage{hyperref} 
\hypersetup{
    colorlinks=true,       
    pagebackref=true, 
    citecolor=blue, 
    linkcolor=black,          
    urlcolor=blue,           
    linktocpage=true, 
    pdfborderstyle={/S/U/W 1}, 
    hyperfootnotes=false 
}

\begin{document}
\title{\huge
Stochastic dynamics of determinantal processes by integration by parts
}
\author{Laurent Decreusefond\thanks{D\'epartement INFRES,
Telecom Paristech, 23 avenue d'Italie, 75013 Paris, France. e-mail: \tt{laurent.decreusefond@mines-telecom.fr}}\and Ian Flint \thanks{Division of Mathematical Sciences,
Nanyang Technological University, SPMS-MAS-05-43, 21 Nanyang Link
Singapore 637371. e-mail: \tt{iflint@ntu.edu.sg}}\and Nicolas Privault\thanks{Division of Mathematical Sciences,
Nanyang Technological University, SPMS-MAS-05-43, 21 Nanyang Link
Singapore 637371. e-mail: \tt{nprivault@ntu.edu.sg}} \and Giovanni
Luca Torrisi\thanks{Istituto per le Applicazioni del Calcolo
"Mauro Picone", CNR, Via dei Taurini 19, 00185 Roma, Italy.
e-mail: \tt{torrisi@iac.rm.cnr.it}} }
\maketitle
\vspace{-0.8cm}
\begin{abstract} 
 We derive an integration by parts formula for 
 functionals of determinantal processes on compact sets,
 completing the arguments of \cite{camilier}.
 This is used to show
 the existence of a configuration-valued diffusion process 
 which is non-colliding 
 and admits the distribution of the determinantal process as
 reversible law. 
 In particular, this approach allows us to build a concrete example of
 the associated diffusion process, providing an illustration of the results of
 \cite{camilier} and \cite{yoo}. 
\end{abstract}


\vspace{0.3cm}


\baselineskip0.7cm

\section{Introduction}\label{sec:introduction}

 Determinantal processes are point processes that exhibit
 repulsion, and were introduced to
 represent the configuration of fermions, cf. \cite{macchi,shirai,tamura}.
 They are known to be connected with the zeros of
 analytic functions (cf. \cite{hough} and references therein) as well as with the theory of random matrices (cf. \cite{anderson}).
To the best of our knowledge, the first use of determinantal processes as models in applications trace back to \cite{LDa}.
More recently, in \cite{flpnw2,MSa,TLa,DFMVa}, different authors have used determinantal processes to model phenomena arising in telecommunication networks.
 \\

 \noindent
 The Markov process associated to the
 Ornstein-Uhlenbeck operator on the Poisson space
 has been constructed in 
 \cite{surgailis}. 
 In \cite{akr}, using the Dirichlet forms theory
 (see \cite{fukushima,rockner}), the diffusion 
 whose symmetrizing measure is the law of a Poisson process
 over a Riemannian manifold is constructed. 
 This result has been extended to Gibbs processes on $\R^d$ in
 \cite{akr2}.
 Using such methods, the Dirichlet form and diffusion process 
 associated to determinantal processes
 has been constructed in \cite{yoo}. 
 \\

 \noindent
 In this paper, by completing the arguments of \cite{camilier}
 we prove an integration by parts formula for 
 functionals of determinantal processes on compact subsets of $\R^d$, 
 and we recover the closability of the associated
 Dirichlet form. 
 This provides a novel proof of the existence of interacting diffusion
 processes properly associated to determinantal processes.
 In addition, our approach based on integration by parts exhibits 
 the generator of the diffusion process,
 and allows in turn to provide an explicit example of
 a diffusion process satisfying our hypotheses.
\\

\noindent
 As a preliminary step, we derive an integration by parts formula
 for functionals of a determinantal process on a compact set $D\subset\R^d$, by completing
 the result established in \cite{camilier}.
 In comparison with \cite{camilier},
 the integration by parts formula on compact sets is extended
 to closed gradient and divergence operators
 by the use of a different set of test functionals,
 cf. \eqref{F}
 and Theorem~\ref{thm:integrationbyparts}.
 Our construction of the diffusion processes
 follows the lines of \cite{akr}, and it
 differs from the one of \cite{torrisi} which is based
 on sample-path identities. 
 Our gradient and divergence operators also
 differ from those of \cite{privaulttorrisi},
 which also deals with compact subsets of $\R^d$.
 Nevertheless the integration by parts formula of
 Theorem~\ref{thm:integrationbyparts} can also be applied
 to density estimation and sensitivity analysis
 for functionals of determinantal processes
 along the same lines. 
\\

\noindent
 In Theorem~\ref{thm:dirichlet} we construct
 the Dirichlet form corresponding to a determinantal process on a compact set $D\subset\R^d$.
 In Theorem~\ref{thm:diff1} we show the existence of the
 diffusion properly associated to a determinantal process
 on a compact set $D\subset\R^d$.
 Note that, as in the other constructions (cf. \cite{akr,akr2,yoo}),
 the associated
 diffusion process admits the distribution of the determinantal process
 as a reversible law.  
 We prove the non-collision property of the
 diffusion in Theorem \ref{thm:noncollision}.
Finally, in Section~\ref{sec:examples}, 
we provide an example of a determinantal process satisfying our integration by parts formula, and for which the aforementioned properly associated diffusion process exists.
\\

 \noindent
 Some definitions related to point processes theory and more in particular
 to determinantal processes are recalled in Section~\ref{sec:preliminaries} 
 based on \cite{daley,daley2,hough,moller,vanlieshout}. 
 Some notions from the Dirichlet forms theory are given in Sections \ref{sec:df} and \ref{sec:diff} based on \cite{fukushima,rockner}.
 We also refer the reader to 
 \cite{dunford,simon2}
 for the required background on functional analysis.
\section{Preliminaries}
\label{sec:preliminaries}
\subsubsection*{Locally finite point processes}\label{subsec:locfinitepp}
 Let $S$ be a 
 Polish space, and denote by $\mathcal{B}(S)$ the associated Borel $\sigma$-algebra. For any subset
$B\subset S$, let $\sharp B$ denote the cardinality of $B$,
setting $\sharp B=\infty$ if $B$ is not finite. We denote by
$\mathrm{N}_{lf}$ the set of locally finite point configurations
on $S$:
\[
\text{$\mathrm{N}_{lf}$:=\{$B\subset S$\,\,\,:\,\,\,$\sharp(B\cap
D)<\infty$, for any compact $D\subset S$\}.}
\]
We identify locally finite configurations with $\mathbb N$-valued simple Radon measures, equip $\mathrm{N}_{lf}$ with the vague topology (see Appendix~2 in \cite{daley}), and we denote the corresponding Borel $\sigma$-algebra by $\mathcal{N}_{lf}$.
We recall that a non-negative simple Radon measure is a Radon measure which is less than or equal to $1$ on singletons.
We define similarly $\mathrm{N}_{f}$ the set of finite point configurations
on $S$:
\[
\text{$\mathrm{N}_{f}$:=\{$B\subset S$\,\,\,:\,\,\,$\sharp B<\infty$\},}
\]
 which is naturally equipped with the trace $\sigma$-algebra $\mathcal{N}_{f} = \mathcal{N}_{lf} |_{\mathrm{N}_{f}}$. For any measurable set $B\subset S$, let $\mathrm{N}_{f}^{B}$
be the space of finite configurations on $B$, and $\mathcal{N}_{f}^B$ the associated trace-$\sigma$-algebra.
\\

\noindent
 By a locally finite and simple point process $\bold X$ on $S$ we mean a
 measurable mapping defined on some probability space
$(\Omega,\mathcal{F},P)$ taking values on
$(\mathrm{N}_{lf},\mathcal{N}_{lf})$. 
We denote by $\bold X(B)$ the number
of points of $\bold X$ in a measurable set $B\subset S$, i.e. $\bold X(B) := \sharp(\bold X \cap B)$, and by
$$
 \bold{X}^{B} =
 \bold{X} \cap B
 =
 \{ X_1,\ldots,X_{\bold X(B)} \}
$$
 the restriction to $B$ of the point process
 $\bold{X}\equiv\{X_n\}_{1\leq n\leq \bold X(S)}$. 
 In the following, we shall denote by $\PP$ the law of $\bold X$ and by $\PP_B$ the law of ${\bold X}^B$.
\\

\noindent
The correlation functions of $\bold X$, with respect to ({\em w.r.t.}) a given
Radon measure $\nu$ on $(S,\mathcal{B}(S))$, are (if they exist) symmetric measurable functions
 $\rho_n:S^n \longrightarrow \R_+$ such that
\[
\mathbb{E}\left[\prod_{i=1}^{n}\bold X(B_i)\right]=\int_{B_1\times\ldots\times
B_n}\rho_n(x_1,\dots,x_n)\,\nu(\mathrm{d}x_1)\cdots\nu(\mathrm{d}x_n),
\]
 for any
 family of mutually disjoint bounded subsets $B_1,\ldots,B_n$ of $S$,
 $n\geq 1$.
 We require in addition that 
 $\rho_n(x_1,\ldots,x_n) =0$ 
 whenever $x_i=x_j$ for some $1 \leq i\neq j \leq n$. When
 $\rho_1$ exists, the measure $\rho_1(x)\,\nu(\mathrm{d}x)$ is known as the intensity measure of $\bold X$.
 \\

\noindent
As in \cite{georgiiyoo}, we define for any Radon measure $\nu$ on $(S,\mathcal{B}(S))$
the $\nu$-sample measure $L^\nu$ on $(\mathrm{N}_{f}, \mathcal{N}_{f})$ by
\begin{equation}
\label{eq:samplemeasure}
\int_{\mathrm{N}_{f}} f(\alpha) \, L^\nu(\mathrm{d}\alpha) := \sum_{n \ge 0} \frac{1}{n!} \int_{S^n} f(\{x_1,\dots,x_n\}) \, \nu(\mathrm{d}x_1)\cdots \nu(\mathrm{d}x_n),
\end{equation}
for any measurable $f:\mathrm{N}_{f} \rightarrow \mathbb{R}_+$.
For any compact subset $D\subset S$, the Janossy densities of $\bold X$ {\em w.r.t.}
$\nu$
are (if they exist) measurable symmetric functions $j^n_{D}:D^n\rightarrow\mathbb{R}$ satisfying, for all measurable $f:\mathrm{N}_{f}^{D} \rightarrow\mathbb{R}_+$,
\begin{equation}
\label{def:janossy}
\mathbb{E}\left[ f(\bold X^D) \right] = \sum_{n \ge 0} \frac{1}{n!} \int_{D^n} f(\{x_1,\dots,x_n\})\,j_{D}^n\left( x_{1},\dots, x_{n}\right)\,\nu(\mathrm{d}x_1)\cdots\nu(\mathrm{d}x_n),
\end{equation}
i.e. defining $j_D(\bold{x}) := j_D^{{\bold x}(D)}(x_1,\dots,x_{{\bold x}(D)})$ for $\bold x=\{x_1,\dots,x_{\bold x(D)}\}\in\mathrm{N}_{f}^{D}$, $j_D$ is the density of $\PP_D$ with respect to $L_D^\nu$ (the restriction to $\mathcal{N}_f^D$ of $L^\nu$), when $\PP_D \ll L_D^\nu$.
\subsubsection*{Kernels and integral operators}
\label{subsec:integraltraceop}

Let $\nu$ be a Radon measure on $(S,\mathcal{B}(S))$. For any compact set $D\subset S$, we denote by $\mathrm{L}^2(D,\nu)$ the
Hilbert space of complex-valued square integrable functions {\em
w.r.t.} the restriction of the Radon measure $\nu$ to $D$,
equipped with the inner product
\[
\langle f,g \rangle_{\mathrm{L}^2(D,\nu)}:=\int_{D}f(x)\overline{g(x)}\,\nu(\mathrm{d}x),\quad\text{$f,g\in
\mathrm{L}^2(D,\nu)$}
\]
where $\overline z$ denotes the complex conjugate of a complex $z\in\C$.
By definition, an integral operator $\mathcal T:\mathrm{L}^2(S,\nu)\to\mathrm{L}^2(S,\nu)$ with kernel $T:S^2\to\C$ is a bounded operator defined by
\[
\mathcal{T}f(x):=\int_{S}T(x,y)f(y)\,\nu(\mathrm{d}y),\quad\text{for
$\nu$-almost all $x\in S$}.
\]
Letting $\mathcal{P}_D$ denote the projection operator
from $\mathrm{L}^2(S,\nu)$ onto $\mathrm{L}^2(D,\nu)$, we set
$\mathcal{T}_D=\mathcal{P}_D \mathcal{T} \,\mathcal{P}_D$ and note that its kernel is $T_D(x,y):={\bold 1}_D(x) T(x,y) {\bold 1}_D(y)$, for $\nu$-almost all $x,y \in S$.
It can be shown that $\mathcal{T}_D$ is a compact operator.
The operator $\mathcal{T}$ is said to be Hermitian or
self-adjoint if
\begin{equation}\label{eq:Hermitian}
T(x,y)=\overline{T(y,x)},\quad\text{for $\nu^{\otimes 2}$-almost
all $(x,y)\in S^2$.}
\end{equation}
Equivalently, this means that the integral operator
$\mathcal{T}_D$ is self-adjoint for any compact set $D\subset S$.
If $\mathcal{T}_D$ is self-adjoint, by the spectral theorem 
we have that $\mathrm{L}^2(D,\nu)$ has
an orthonormal basis $( \varphi_j^D )_{j\geq 1}$ of eigenfunctions
of $\mathcal{T}_D$. The corresponding eigenvalues
$( \lambda_j^D )_{j\geq 1}$ have finite multiplicity (except
possibly the zero eigenvalue) and the only possible accumulation
point of the eigenvalues is the eigenvalue zero.
 Then, the kernel $T_D$ of $\mathcal{T}_D$ can be written as
\begin{equation}
\label{eq:kdecomp}
T_D(x,y) = \sum_{j \ge 1} \lambda_j^D \varphi_j^D(x) \overline{ \varphi_j^D(y) },
\end{equation}
for $\nu$-almost all $x, y \in D$. We say that an operator $\mathcal{T}$ is positive (respectively non-negative) if its spectrum is included in $(0,+\infty)$ (respectively $[0,+\infty)$). For two operators $\mathcal{T}$ and $\mathcal{U}$, we will say that $ \mathcal{T}> \mathcal{U}$ (respectively $ \mathcal{T} \ge \mathcal{U}$) in the operator ordering if $\mathcal{T}-\mathcal{U}$ is a positive operator
(respectively a non-negative operator).
\\

\noindent
 We say that a self-adjoint integral operator $\mathcal{T}_D$, with kernel $T_D$ as in \eqref{eq:kdecomp}, is of
trace class if
\[
\sum_{j\geq 1}|\lambda_j^D|<\infty,
\]
 and we
 define the trace of the operator $\mathcal{T}_D$ as
$$
 \mathrm{Tr\,}(\mathcal{T}_D) := \sum_{j\geq 1}\lambda_j^D.
$$
 If $\mathcal{T}_D$ is of trace class for every compact subset
$D\subset S$, then we say that $\mathcal{T}$ is locally of trace
class. It is easily seen that $\mathcal{T}^n$ is locally of trace class, for all $n \ge 2$, if $\mathcal{T}$ is locally of trace class.
Finally, we define the Fredholm determinant of $\bold{Id}+\mathcal T_D$, when $\| \mathcal{T}_D \|_{op} < 1$, as
\begin{equation}
\label{eq:fredholm}
\mathrm{Det}(\bold{Id}+\mathcal{T}_D) := \mathrm{exp}\left(\sum_{n \ge 1}\frac{(-1)^{n-1}}{n} \mathrm{Tr}(\mathcal{T}_D^n) \right).
\end{equation}
Here, $\bold{Id}$ denotes the identity operator on $\mathrm{L}^2(S,\nu)$ and $\|\cdot\|_{op}$ denotes the operator norm.
\subsubsection*{Determinantal processes on $S$}
\label{sec:determinpp}
Let $\mu$ be a Radon measure on $(S,\mathcal B(S))$.
A locally
finite and simple point process $\bold{X}$ 
on $S$ is said to be a determinantal process with kernel $K$ and reference measure $\mu$ if its
correlation functions {\em w.r.t.} $\mu$ exist and
satisfy
\[
\rho_n(x_1,\ldots,x_n)=\mathrm{det}(K(x_i,x_j))_{1\leq i,j\leq n},
\]
\noindent for any $n\geq 1$ and $\mu$-a.e. $x_1,\ldots,x_n\in S$, where $K(\cdot,\cdot)$ is a measurable function.
Throughout this paper we shall work under the following
hypothesis.
\\

\noindent
\text{\bf (H1)}: {\it
 The integral operator $\mathcal{K}$ on $\mathrm L^2(S,\mu)$ with kernel $K$ is locally of trace class, self-adjoint,
and $0 \le \mathcal{K} < \bold{Id}$.}
\\

The existence and uniqueness in law of a determinantal process with kernel $K$ is
guaranteed under \text{\bf (H1)} by the results in \cite{macchi,shirai,soshnikov}. 
See also Lemma 4.2.6 and Theorem 4.5.5 in
\cite{hough}.
\\

\noindent
 We define the global interaction operator $\mathcal{J} := (\bold{Id} - \mathcal{K})^{-1} \mathcal{K}$.
 As proved in \cite{georgiiyoo}, 
 we also define the local interaction operator
$\mathcal{J}[D]$ as 
$$\mathcal{J}[D]:=(\bold{Id}-\mathcal{K}_D)^{-1}\mathcal{K}_D,
$$
 and we emphasize that
 $\mathcal{J}[D]$ is not the projection of $\mathcal J$ onto $\mathrm L^2(D,\mu)$.
 We refer the reader to \cite{georgiiyoo}
 for the following properties of $\mathcal{J}[D]$.
 First, $\mathcal{J}[D]$ is a self-adjoint integral operator and, letting
$J[D]$ denote its kernel,
as a consequence of \eqref{eq:kdecomp} we have
\begin{equation}
\label{eq:jdecomp}
J[D](x,y) = \sum_{j \ge1} \frac{\lambda_j^D}{1-\lambda_j^D}\, \varphi_j^D(x) \overline{ \varphi_j^D(y) },
\end{equation}
for $\mu$-almost all $x, y \in D$. Second, 
$\mathcal{J}[D]$ is a trace class operator.
Third, denoting by
$\mathrm{det}\,J[D](\{x_1,\ldots,x_n\}):=\mathrm{det}\,(J[D](x_i,x_j))_{1\leq i,j\leq
n}$, the function
\[
(x_1,\ldots,x_n)\mapsto
\mathrm{det}\,J[D](\{x_1,\ldots,x_n\})
\]
is $\mu^{\otimes n}$-a.e. non-negative 
and symmetric in
$x_1,\ldots,x_n$, 
and we set
$\mathrm{det}\,J[D](\{x_1,\ldots,x_n\}):=\mathrm{det}\,J[D](x_1,\ldots,x_n)$.
The local interaction operator is related to the Janossy densities of
a determinantal process by the following proposition.

\begin{proposition}[\cite{shirai}]
\label{prop:janossydensityDPP}
Under \text{\bf (H1)}, for any compact $D \subset S$,
 the Janossy densities $j_D^n(\bold{x})$ of $\bold X$ are given by
\begin{equation}
\label{eq:janossy}
j_D^n(x_1,\dots,x_n) = \mathrm{Det } (\bold{Id} - \mathcal{K}_D)\, \mathrm{det}\,J[D](x_1,\dots,x_n),
 \qquad
x_1,\ldots,x_n\in D,
 \quad n \geq 1.
\end{equation}
 Moreover, $P(\bold X^D=\emptyset)$ is given by
$j_D^0 (\emptyset) = \mathrm{Det } (\bold{Id} - \mathcal{K}_D)$.
\end{proposition}
\section{Differential calculus and integration by parts}
\label{sec:diffintpartconfigspace}
 In this section we derive an
 integration by parts formula for functionals of a determinantal point process, and we extend it by closability.
 Hereafter we assume that $S$ is a domain of $\R^d$
 equipped with the Euclidean distance,
$\mu$ is a Radon measure on $(S,\mathcal B(S))$
 and $D\subset S$ is a fixed compact set. We denote by
$\|\cdot\|$ the Euclidean norm on $\R^d$, by $x\cdot y$ the usual inner product of $x,y\in\R^d$, and by $x^{(i)}$ the
$i$-th component of $x\in\R^d$, $i=1,\ldots, d$.
\subsubsection*{Differential calculus}
\label{sec:diffconfigspace}
 We denote by
 ${\cal C}^\infty (D ,\R^d)$ the set of all
 ${\cal C}^\infty$-vector fields $v:D \longrightarrow \R^d$
and by ${\cal C}^\infty(S^k )$
the set of all
${\cal C}^\infty$-functions on $S^k$.
\begin{definition}
\label{def:testfunctions}
A function $F:\mathrm N_f^D\to\R$ is said to be in ${\cal
S}_D$
if
\begin{equation}
\label{F}
 F (\bold{x} )
 =
 f_0
 {\bf 1}_{\{ \bold x ( D ) = 0 \}}
 +
 \sum_{k=1}^n
 {\bf 1}_{\{ \bold x ( D ) = k \}}
 f_k( x_1,\ldots,x_k )
,\quad\text{ for }\mathbf x=\{x_1,\dots,x_{\mathbf x(D)}\}\in\mathrm N_f^D,
\end{equation}
where $f_0\in \R$ is
a constant, $n\geq 1$ is an integer and, for any $k=1,\ldots,n$, $f_k \in
{\cal C}^\infty(D^k)$ is a symmetric function.
\end{definition}
\noindent
The set of test functions ${\cal
S}_D$ is dense in $\mathrm{L}_D^2:=\mathrm{L}^2(\mathrm{N}_f^D,\PP_D)$ (indeed, it contains the space $\tilde{\mathcal S}_D$ defined in Definition~\ref{def:testfunctions2} below which is dense in $\mathrm{L}_D^2$, see e.g. \cite{rockner} p.54).
 
\noindent
 The gradient of $F\in\mathcal{S}_D$ as in \eqref{F}
 is
 defined by
\begin{equation}
\label{eq:nablax}
 \nabla_x^{\mathrm{N}_{lf}} F (\bold{x} ) 
 :=
 \sum_{k=1}^n
 {\bf 1}_{\{ \bold x ( D ) = k \}}
 \sum_{i=1}^k
 {\bf 1}_{\{ x_i \}} (x)
 \nabla_{x_i}
 f_k( x_1,\ldots,x_k )
, \qquad
 x \in D,\ \mathbf x\in\mathrm N_f^D,
\end{equation}
 where $\nabla_{y}$ denotes the usual gradient on $\R^d$
 with respect to $y$.
 For $v\in\mathcal{C}^\infty(D,\R^d)$, we also let
\begin{equation}\label{eq:gradF}
 \nabla_v^{\mathrm{N}_{lf}}
 F (\bold{x})
 :=
 \sum_{k=1}^{\bold x (D)}
 \nabla_{x_k}^{\mathrm{N}_{lf}}F (\bold{x} )
 \cdot
 v(x_k)
 =
 \sum_{k=1}^n
 {\bf 1}_{\{ \bold x ( D ) = k \}}
 \sum_{i=1}^k
 \nabla_{x_i}
 f_k( x_1,\ldots,x_k )
 \cdot
 v(x_i)
,
\end{equation}
 where we recall that the symbol $\cdot$ denotes the inner product on $\R^d$.
\subsubsection*{Quasi-invariance}
 Next we recall some results from \cite{camilier},
 with some complements that make the proofs more precise.
 Let $\mathrm{Diff} _{0}(D)$ be the set of all diffeomorphisms
 from $D$ into itself.
 For $\phi \in \mathrm{Diff}_0(D)$ and a Radon measure $\nu$ on $D$, $\nu_{\phi}$ denotes
the image measure of $\nu$ by $\phi$.
For such a $\phi$, we define the
map 
\begin{align*}
\mathcal{I}_\phi \, :\, \mathrm{L}^2(D, \nu_{\phi})& \longrightarrow \mathrm{L}^2(D,\nu),\\*
f &\longmapsto f\circ \phi,
\end{align*}
whose inverse is given by $\mathcal{I}_\phi^{-1}=\mathcal I_{\phi^{-1}}$. Note that $\mathcal{I}_\phi $ is an isometry. 
Given an operator $\mathcal T$ on $\mathrm L^2(D,\nu)$, we define the operator on $\mathrm L^2(D,\nu_\phi)$
$$
\mathcal{T}^{\phi}=\mathcal{I}_\phi^{-1}\,\mathcal{T}\,\mathcal{I}_\phi. 
$$
Lastly, for any $\bold x = \{{x}_n\}_{1\leq n\leq \bold x(S)} \in \mathrm{N}_{lf}$, we denote by $\Phi$
the map:
\begin{align*}
  \Phi\, :\,  \mathrm{N}_{lf} & \longrightarrow  \mathrm{N}_{lf},\\*
 \{{x}_n\}_{1\leq n\leq \bold x(S)}&\longmapsto \{\phi({x}_n)\}_{1\leq n\leq \bold x(S)} .
\end{align*}

\noindent
The following lemma is proved in \cite{camilier} and \cite{torrisi4}.
\begin{lemma}
  \label{lem:operators}
  Assume \text{\bf (H1)},
  and take $\phi\in\mathrm{Diff}_0(D)$. The following properties hold.
  \begin{enumerate}[a)]
  \item $\mathcal{K}^{\phi}_{D}$ and $\mathcal{J}[D]^{\phi}$ are integral
  operators on $\mathrm{L}^2(D,\mu_\phi)$ with kernels given respectively by
$K_D^\phi(x,y)=
  K\left(\phi^{-1}(x),\phi^{-1}(y)\right)
$
and $J[D]^\phi(x,y)=
  J[D] \left(\phi^{-1}(x),\phi^{-1}(y)\right)$.
  \item $\mathcal{K}^{\phi}_{D}$ is of trace class and
    $\mathrm{Tr}(\mathcal{K}^{\phi}_{D})=\mathrm{Tr}(\mathcal{K}_D)$.
  \item $\mathrm{Det}(\bold{Id}-
    \mathcal{K}_{D}^{\phi})=\mathrm{Det}(\bold{Id}-\mathcal{K}_{D})$. This translates into the fact that $P(\bold X (D) = 0) = P( \Phi(\bold X) (D) = 0) $ which is expected since $\phi$ is a diffeomorphism.
  \item $\mathcal{J}[D]^{\phi}=\mathcal J^\phi[D]:=(\bold{Id} - \mathcal{K}^{\phi}_D)^{-1} \mathcal{K}^{\phi}_D$ is the local interaction operator associated with $\mathcal{K}^{\phi}$.
  \end{enumerate}
\end{lemma}
The following mapping theorem holds, see Theorem~7 in \cite{camilier}.
\begin{lemma}
\label{lem:phidet}
Assume \text{\bf (H1)}, and let $\phi\in\mathrm{Diff}_0(D)$.
Then, $\Phi(\bold X^D)$ is a determinantal process with integral operator $\mathcal{K}^\phi_D$ and reference measure $\mu_\phi$.
\end{lemma}

\noindent
To prove the quasi-invariance of the determinantal measure restricted to a compact set $D \subset S$ with respect to the group of diffeomorphisms on $D$, we state one last result.
\begin{lemma}
\label{lem:detj}
Under \text{\bf (H1)}, we have $\mathrm{det}\,J[D](\bold{x})>0$, for $\PP_D$-a.e. $\bold{x}\in\mathrm{N}_{f}^D$. However, it does not in general hold that $\mathrm{det}\,J[D](\bold{x})>0$, for $L_D^\mu$-a.e. $\bold{x}\in\mathrm{N}_{f}^D$.
\end{lemma}
\begin{Proof}
Recall that under \text{\bf (H1)}, we have $\PP_D \ll L_D^\mu$ and
$$
j_D(\bold{x}) = \frac{\mathrm{d}\PP_D}{\mathrm{d}L_D^\mu}(\bold{x}) = \mathrm{Det } (\bold{Id} - \mathcal{K}_D)\, \mathrm{det}\,J[D](\bold{x}),
$$
for $\bold x \in \mathrm{N}_{f}^D$. Since $j_D$ is a density, we clearly have $j_D(\bold x) > 0$, for $\PP_D$-a.e. $\bold x \in \mathrm{N}_{f}^D$. Hence, since $\| \mathcal{K}_D \| < 1$, we have  $\mathrm{det}\,J[D](\bold{x})>0$, for $\PP_D$-a.e. $\bold x \in \mathrm{N}_{f}^D$. As to the concluding part of the lemma, we notice that in general, one does not have that $\PP_D$ is equivalent to  $L_D^\mu$. Indeed, consider for example the case where the rank of $\mathcal{K_D}$ is less than or equal to $N\ge1$. Then, $j_D^{N+1}(x_1,\dots,x_{N+1}) = 0$, for $\mu^{\otimes (N+1)}$-a.e. $(x_1,\dots,x_{N+1})\in D^{N+1}$ (since $\bold{X}^D$ has less than $N+1$ points almost surely, see \cite{soshnikov} for details). It suffices to define the set
$$
A := \{B\subset D\,\,\,:\,\,\,\sharp B = N+1\},
$$
which verifies $\PP(A) = 0$ but $L_D^\mu(A) = \frac{1}{(N+1)!} \mu(D)$.
\end{Proof}
\begin{remark}
If we assume that,
for any $n\geq 1$, the function
$$
(x_1,\ldots,x_n)\longmapsto\mathrm{det}\,J[D](x_1,\ldots,x_n)
$$
is strictly positive $\mu^{\otimes n}$-a.e. on $D^n$, then we have that $\PP_D$ and $L_D^\mu$ are equivalent, and it follows that $\mathrm{det}\,J[D](\bold{x})>0$, for $L_D^\mu$-a.a. $\bold{x}\in\mathrm{N}_{f}^D$.
\end{remark}

\noindent
 The next Proposition~\ref{djklasda1} is similar
 to its analog in \cite{camilier}, however
 the proof given there implicitly uses the fact that
 $\mathrm{det}\,J[D](\bold{x})>0$, for $L_D^\mu$-a.a.
 $\bold{x}\in\mathrm{N}_{f}^D$, which has been shown to be false in general.
In order to prove Proposition~\ref{djklasda1}, we assume the following technical condition.
\\

\noindent
\text{$\bold{(H2):}$}\,\,{\it The Radon measure $\mu$
 is absolutely continuous {\em w.r.t.} the Lebesgue measure $\ell$
 on $D$, with Radon-Nikodym derivative
 $\rho = \frac{\mathrm{d} \mu}{\mathrm{d} \ell}$ which is strictly positive and continuously differentiable on $D$.
}
\\

\noindent
 Under \text{$\bold{(H2)}$}, for any $\phi \in \mathrm{Diff}_0(D),$ $\mu_\phi$ is absolutely
continuous with respect to $\mu$ with density given by
\begin{equation}
\label{eq:densitymuphi}
  p^{\mu}_{\phi}(x)=\frac{\mathrm{d}\mu_{\phi}(x)}{\mathrm{d}\mu(x)}
  =\frac{\rho(\phi^{-1}(x))}{\rho(x)}\mathrm{Jac}(\phi^{-1})(x),
\end{equation}
where $\mathrm{Jac}(\phi^{-1})(x)$ is the Jacobian of $\phi^{-1}$ at a point $x\in S$. We draw attention to the fact that it is indeed $\mathrm{Jac}(\phi^{-1})(x)$ that appears in (\ref{eq:densitymuphi}),
which differs from equation (2.11) of \cite{akr}.
We are now in a position to state and prove the main result of this section.
\begin{proposition}
\label{djklasda1}
Assume \text{\bf (H1)}, \text{\bf (H2)}, and let $\phi\in\mathrm{Diff}_0(D)$, for any measurable non-negative $f$ on $D$ we have
  \begin{multline*}
  \label{eq:quasiinvariance}
    \mathbb{E}\left[\exp\left(-\sum_{k=1}^{\bold X(D)} f \circ \phi (X_k) \right)\right]=\\
    \mathbb{E}\left[\exp\left(-\sum_{k=1}^{\bold X(D)} f(X_k) +\sum_{k=1}^{\bold X(D)}
        \mathrm{ln}(p^{\mu}_{\phi} (X_k))\right)\frac{\mathrm{det}\,
        J^{\phi}[D]({\bold X}^D)}{\mathrm{det}\, J[D]({\bold X}^D)}\right].
  \end{multline*}
\end{proposition}
\begin{Proof}
For any measurable non-negative $f$ on $D$, we have
\begin{multline*}
 \shoveleft{ \mathbb{E}\left[\exp\left(-\sum_{k=1}^{\bold X(D)} f(X_k) +\sum_{k=1}^{\bold X(D)}
        \mathrm{ln}(p^{\mu}_{\phi} (X_k))\right)\frac{\mathrm{det}\,
        J^{\phi}[D]({\bold X}^D)}{\mathrm{det}\, J[D]({\bold X}^D)}\right]} \\
        \shoveleft{= \sum_{n \ge 0} \frac{1}{n!} \int_{D^n}  \mathrm{e}^{-\sum_{k=1}^n f(x_k) }\prod_{k=1}^n p^{\mu}_{\phi} (x_k)} \\
        \shoveright{\times \frac{\mathrm{det}\,
        J^{\phi}[D](x_1,\dots,x_n)}{\mathrm{det}\, J[D](x_1,\dots,x_n)}\,j_{D}\left( x_{1},\dots, x_{n}\right)
        \,\mu(\mathrm{d}x_1)\cdots\mu(\mathrm{d}x_n)}\\
         \shoveleft{= \sum_{n \ge 0} \frac{1}{n!} \int_{D^n}  \mathrm{e}^{-\sum_{k=1}^n f(x_k) }\prod_{k=1}^n p^{\mu}_{\phi} (x_k)\, \mathrm{det}\,
        J^{\phi}[D](x_1,\dots,x_n)\,\mathrm{Det } (\bold{Id} - \mathcal{K}_D)\,\mu(\mathrm{d}x_1)\cdots\mu(\mathrm{d}x_n)}\\
         \shoveleft{= \sum_{n \ge 0} \frac{1}{n!} \int_{D^n}  \mathrm{e}^{-\sum_{k=1}^n f(x_k) } \, \mathrm{Det } (\bold{Id} - \mathcal{K}_D^\phi) \, \mathrm{det}\,
        J^{\phi}[D](x_1,\dots,x_n)\,\mu_\phi(\mathrm{d}x_1)\cdots\mu_\phi(\mathrm{d}x_n)},\\
\end{multline*}
where we have used $(\ref{eq:janossy})$, $(\ref{eq:densitymuphi})$ and Lemma \ref{lem:operators}, $c)$.
Then, we conclude by Lemma \ref{lem:phidet}. Indeed, $\mathrm{Det } (\bold{Id} - \mathcal{K}_D^\phi) \, \mathrm{det}\, J^{\phi}[D](x_1,\dots,x_n)$ is the Janossy density of $\Phi(\bold X^D)$ with respect to $\mu_\phi$ (see Lemma~\ref{lem:operators}, $e)$).
\end{Proof}
\subsubsection*{Integration by parts and closability}
\label{sec:intbyparts}
 We close this section with the statement and proof
 of the integration by parts formula
 for determinantal processes which is based on the closed gradient
 and divergence operators, cf. Theorem~\ref{thm:integrationbyparts} below.
The integration by parts formula is proved
 on the set of test functionals $\mathcal{S}_D$
 introduced in Definition~\ref{def:testfunctions}, extending and making more precise
 the argument and proof of Theorem 10 page 289 of \cite{camilier}.
\\

\noindent \noindent \text{$\bold{(H3):}$}\,\,{\it 
The function
$
(x_1,\ldots,x_n)\longmapsto\mathrm{det}\,J[D](x_1,\ldots,x_n)
$
 is continuously differentiable on $D^n$.}
\\

\noindent
 Assuming that
\text{$\bold{(H1)}$} and \text{$\bold{(H3)}$}
 hold, the
 potential energy is the function $U:\mathrm{N}_{f}^D
\longrightarrow \R$ defined by
$$
U[D](\bold{x}) :=
-\log\mathrm{det}\,J[D](\bold{x}).
$$
We insist that since $\mathrm{det}\,J[D](\bold{x}) > 0$ for $\PP_D$-a.e.
$\bold{x}\in \mathrm{N}_{f}^D$, $U$ is well defined for $\PP_D$-a.e. $\bold{x}\in \mathrm{N}_{f}^D$.
 For any $v\in {\cal C}^\infty (D,\R^d )$ and $\bold x=\{x_1,\dots,x_{\bold x(D)}\}\in\mathrm N_{f}^D$, we set
\begin{eqnarray}
\nonumber
 \nabla_v^{\mathrm{N}_{lf}}
 U[D] (\bold{x} )
 & := &
 -
 \sum_{k=1}^\infty
 {\bf 1}_{\{ \bold x ( D ) = k \}}
 \sum_{i=1}^k
\frac{ \nabla_{x_i}
 \mathrm{det}\,J[D]( x_1,\ldots,x_k )
}{\mathrm{det}\,J[D]( x_1,\ldots,x_k ) }
 \cdot
 v(x_i)
\\
\label{U}
 & = &
 \sum_{k=1}^\infty
 {\bf 1}_{\{ \bold x ( D ) = k \}}\sum_{i=1}^{k}
 U_{i,k}( x_1,\ldots,x_k )\cdot v(x_i)
.
\end{eqnarray}
\noindent
 Under Conditions $\bold{(H1)}$ and $\bold{(H2)}$
 we define the vector field
\begin{equation*}
\beta^{\mu}(x):=
 \frac{\nabla\rho(x)}{\rho(x)},
\end{equation*}
as well as 
\[
B_v^{\mu}(\bold{x}):=\sum_{k=1}^{\bold x (D)}(-\beta^\mu(x_k)\cdot
v(x_k)+\div v(x_k)),
\]
 where $\div$ denotes the adjoint of the gradient $\nabla$ on $D$.
\begin{lemma}
\label{le:integrationbyparts}
 Assume $\bold{(H1)}$,
 $\bold{(H2)}$ and $\bold{(H3)}$. 
 Then, for any
$F,G\in {\cal S}_D$ and $v\in {\cal C}^\infty (D,\R^d )$, we have
\begin{equation}
\label{ibpl}
\mathbb{E}[G(\bold{X}^{D})\nabla_{v}^{\mathrm{N}_{lf}}F(\bold{X}^{D})]
=\mathbb{E}[F(\bold{X}^{D}) \nabla_v^{\mathrm{N}_{lf}*} G (\bold{X}^{D})],
\end{equation}
 where
\[
\nabla_v^{\mathrm{N}_{lf}*} G (\bold{x}) :=-\nabla_{v}^{\mathrm{N}_{lf}}G(\bold{x})+G(\bold{x})
\left( - B_v^\mu(\bold{x})+\nabla_{v}^{\mathrm{N}_{lf}}U[D](\bold{x})\right),\qquad\bold x\in\mathrm N_f^D.
\]
\end{lemma}
\begin{Proof}
 For $v\in {\cal C}^\infty (D,\R^d )$,
 consider the flow $\phi_t^v:D \longrightarrow D$, $t\in\R$, where for a fixed $x\in D$,
 the curve $t\mapsto\phi_t^v(x)$ is defined as the solution to the
 Cauchy problem
\[
 \frac{\mathrm{d}}{\mathrm{d}t}\phi_t^v(x)=v(\phi_t^v(x)),\quad\phi_0^v(x)=x.
\]
\noindent 
 We define the mapping
 $\Phi_t^v:\mathrm{N}_{f}^D \longrightarrow \mathrm{N}_{f}^D$ by
$
\Phi_t^v(\bold{x}):=\{\phi_t^v(x):\,\,x\in\bold{x}\}.
$
Following \cite{akr}, for a function $R:\mathrm N_f^D\to\R$, we define the gradient
$\nabla_v^{\mathrm{N}_{lf}} R(\bold{x})$ as the directional
derivative along $v$ i.e. for $\bold x\in\mathrm N_f^D$,
\[
\nabla_v^{\mathrm{N}_{lf}}
R(\bold{x}):=\frac{\mathrm{d}}{\mathrm{d}t}R(\Phi_t^v(\bold{x}))\Big|_
{t=0},
\]
provided that the derivative exists. It is easy to check that formulas
\eqref{eq:gradF} and \eqref{U} are consistent with this
definition. Note that by \eqref{eq:densitymuphi}, the image measure $\mu_{\phi_{t}^v}$ is
absolutely continuous with respect to $\mu$ on $D$, with
Radon-Nikodym derivative
\[
\frac{\rho(\phi_{-t}^v(x))}{\rho(x)}\,\mathrm{Jac}\left({\phi_{-t}^v}\right)(x).
\]
Note also that
\begin{equation}\label{eq:jacob}
\mathrm{Jac}\left({\phi_t^v}\right)(x)=\exp\left(-\int_{0}^{t}\mathrm{div}\,v(\phi_{z}^v(x))\,\mathrm{d}z\right), 
\end{equation}
and therefore
\begin{align}
\frac{\mathrm{d}}{\mathrm{d}t}\left(\frac{\rho(\phi_{-t}^v(x))}{\rho(x)}\,\mathrm{Jac}\left({\phi_{-t}^v}\right)(x)\right)
=&-\exp\left(-\int_{0}^{t}\mathrm{div}\,v(\phi_{z}^v(x))\,\mathrm{d}z\right) \Big[ \frac{\nabla\rho(\phi_{-t}^v(x))}{\rho(x)}\cdot
v(\phi_{-t}^v(x))\nonumber\\
&+\frac{\rho(\phi_{-t}^v(x))}{\rho(x)}\,\mathrm{div}\,v(\phi_{t}^v(x)) \Big].\label{eq:derjac}
\end{align}
Using Proposition~\ref{djklasda1}, for any $t\in\R$ and
$F,G\in\mathcal{S}_D$, we have
\begin{align}
\mathbb{E}[F(\Phi_{t}^v(\bold{X}^D))&G(\bold{X}^D)]=\mathbb{E}\left[F(\bold{X}^D)G(\Phi_{-t}^v(\bold{X}^D))\vphantom{\left(\prod_{k=1}^{\bold X (D)}
\frac{\rho(\phi_{-t}^v(X_k))}{\rho(X_k)}\,\mathrm{Jac}\left({\phi_{-t}^v}\right)(X_k)\right)}\right.\label{jkldklddd}
\\
&\left.\times\left(\prod_{k=1}^{\bold X (D)}
\frac{\rho(\phi_{-t}^v(X_k))}{\rho(X_k)}\,\mathrm{Jac}\left({\phi_{-t}^v}\right)(X_k)\right)
\frac{\mathrm{det}\,J[D](
\phi_{-t}^v(X_1),\ldots,\phi_{-t}^v(X_{\bold X (D)})
)}{\mathrm{det}\,J[D]( X_1,\ldots,X_{\bold X (D)} )}\right].\nonumber
\end{align}
 We now differentiate this relation with respect to $t$, and exchange 
 $\mathrm{d}/\mathrm{d}t$ with 
 $\mathbb E$. 
 This exchange will be justified later on after
 \eqref{eq:derminuslog}. 
Writing, for ease of notation, $\mathrm{Jac}^{\phi_t^v}:=\mathrm{Jac}\left({\phi_t^v}\right)$, we have
\begin{align}
\label{eq:eqthm36(1)}
&\mathbb{E}\left[G(\bold{X}^D)\frac{\mathrm{d}}{\mathrm{d}t}F(\Phi_{t}^v(\bold{X}^D))\right]\\
&=\mathbb{E}\left[\left(\frac{\mathrm{d}}{\mathrm{d}t}G(\Phi_{-t}^v(\bold{X}^D))\right)
F(\bold{X}^D)\left(\prod_{k=1}^{\bold X (D)}
\frac{\rho(\phi_{-t}^v(X_k))}{\rho(X_k)}\,\mathrm{Jac}^{\phi_{-t}^v}(X_k)\right)
\right.\nonumber\\
&\qquad\qquad\qquad\qquad\qquad\qquad\qquad\qquad\qquad\qquad\left.\vphantom{\left(\prod_{k=1}^{\bold X (D)}
\frac{\rho(\phi_{-t}^v(X_k))}{\rho(X_k)}\,\mathrm{Jac}^{\phi_{-t}^v}(X_k)\right)}
\times\frac{\mathrm{det}\,J[D](
\phi_{-t}^v(X_1),\ldots,\phi_{-t}^v(X_{\bold X (D)})
)}{\mathrm{det}\,J[D]( X_1,\ldots,X_{\bold X (D)} )}\right]\label{eq:eqthm36(2)}\\
&+\mathbb{E}\left[\left(\frac{\mathrm{d}}{\mathrm{d}t}\prod_{k=1}^{\bold X (D)}
\frac{\rho(\phi_{-t}^v(X_k))}{\rho(X_k)}\,\mathrm{Jac}^{\phi_{-t}^v}(X_k)\right)
F(\bold{X}^D)G(\Phi_{-t}^v(\bold{X}^D))
\right.\nonumber\\
&\qquad\qquad\qquad\qquad\qquad\qquad\qquad\qquad\qquad\qquad\left.\vphantom{\left(\frac{\mathrm{d}}{\mathrm{d}t}\prod_{k=1}^{\bold X (D)}
\frac{\rho(\phi_{-t}^v(X_k))}{\rho(X_k)}\,\mathrm{Jac}^{\phi_{-t}^v}(X_k)\right)}
\times\frac{\mathrm{det}\,J[D](
\phi_{-t}^v(X_1),\ldots,\phi_{-t}^v(X_{\bold X (D)})
)}{\mathrm{det}\,J[D]( X_1,\ldots,X_{\bold X (D)} )}\right]\label{eq:eqthm36(3)}\\
&+\mathbb{E}\left[\left(\frac{\mathrm{d}}{\mathrm{d}t}\frac{\mathrm{det}\,J[D](
\phi_{-t}^v(X_1),\ldots,\phi_{-t}^v(X_{\bold X (D)})
)}{\mathrm{det}\,J[D]( X_1,\ldots,X_{\bold X (D)} )}\right)
F(\bold{X}^D)G(\Phi_{-t}^v(\bold{X}^D))\right.\nonumber\\
&\qquad\qquad\qquad\qquad\qquad\qquad\qquad\qquad\qquad\qquad\left.\vphantom{\left(\frac{\mathrm{d}}{\mathrm{d}t}\frac{\mathrm{det}\,J[D](
\phi_{-t}^v(X_1),\ldots,\phi_{-t}^v(X_{\bold X (D)})
)}{\mathrm{det}\,J[D]( X_1,\ldots,X_{\bold X (D)} )}\right)}
\times\prod_{k=1}^{\bold X (D)}
\frac{\rho(\phi_{-t}^v(X_k))}{\rho(X_k)}\,\mathrm{Jac}^{\phi_t^v}(X_k)
\right]\label{eq:eqthm36(4)}.
\end{align}
The claimed integration by parts formula follows by evaluating the
above relation at $t=0$. In particular, we use \eqref{eq:derjac}
to evaluate \eqref{eq:eqthm36(3)}, and we use the relation
\begin{align}
\label{eq:derminuslog}
&\frac{\mathrm{d}}{\mathrm{d}t}\frac{\mathrm{det}\,J[D](
\phi_{-t}^v(X_1),\ldots,\phi_{-t}^v(X_{\bold X (D)})
)}{\mathrm{det}\,J[D]( X_1,\ldots,X_{\bold X (D)} )}\nonumber\\
&=-\sum_{i=1}^{\bold X (D)} \frac{ \nabla_{X_i}
 \mathrm{det}\,J[D](\phi_{-t}^{v}(X_1),\ldots,\phi_{-t}^{v}(X_{\bold X (D)}) )
}{\mathrm{det}\,J[D]( X_1,\ldots,X_{\bold X (D)} ) }
 \cdot
 v(\phi_{-t}^v(X_i))
\end{align}
to evaluate \eqref{eq:eqthm36(4)}.
 Using the definition of $\mathcal{S}_D$, one checks that for $\PP_D$-a.e. $\bold x\in\mathrm N_f^D$,
\[
t\mapsto G(\bold{x})\frac{\mathrm{d}}{\mathrm{d}t}F(\Phi_{t}^v(\bold{x})),
\]
is uniformly bounded by a positive constant in an neighborhood of zero.
 By the assumptions $\bold{(H2)}$ and $\bold{(H3)}$ and the form
 \eqref{F} of the functionals in $\mathcal{S}_D$, one may easily check that
\eqref{eq:eqthm36(2)}, \eqref{eq:eqthm36(3)} and \eqref{eq:eqthm36(4)} can be uniformly bounded in an neighborhood of zero by
$\PP_D$-integrable functions.
This justifies the exchange of derivative and
 expectation in \eqref{jkldklddd}.
 We check this fact only for \eqref{eq:eqthm36(4)}.
Take
$$
F (\bold{x})
 =
 f_0
 {\bf 1}_{\{ \bold x ( D ) = 0 \}}
 +
 \sum_{k=1}^n
 {\bf 1}_{\{ \bold x ( D ) = k \}}
 f_k( x_1,\ldots,x_k ),\qquad\bold x=\{x_1,\dots,x_{\bold x(D)}\}\in\mathrm N_f^D,
$$
 of the form \eqref{F}. By \eqref{eq:derminuslog} we easily see
 that, up to a positive constant, the modulus of the r.v.
\[
\left(\frac{\mathrm{d}}{\mathrm{d}t}\frac{\mathrm{det}\,J[D](
\phi_{-t}^v(X_1),\ldots,\phi_{-t}^v(X_{\bold X (D)})
)}{\mathrm{det}\,J[D]( X_1,\ldots,X_{\bold X (D)} )}\right)
F(\bold{X}^D)G(\Phi_{-t}^v(\bold{X}^D))\prod_{k=1}^{\bold X (D)}
\frac{\rho(\phi_{-t}^v(X_k))}{\rho(X_k)}\,\mathrm{Jac}^{\phi_{-t}^v}(X_k)
\]
is bounded above by
\begin{multline}
\sum_{k=1}^n
 {\bf 1}_{\{ \bold X ( D ) = k \}}\left(\prod_{i=1}^{k}
\frac{\rho(\phi_{-t}^v(X_i))}{\rho(X_i)}\,\mathrm{Jac}^{\phi_{-t}^v}(X_i)\right)
\\
\times\sum_{i=1}^{k}\Big|\frac{
\nabla_{X_i}
 \mathrm{det}\,J[D](\phi_{-t}^{v}(X_1),\ldots,\phi_{-t}^{v}(X_{k}) )
}{\mathrm{det}\,J[D]( X_1,\ldots,X_{k} ) }
 \cdot
 v(\phi_{-t}^v(X_i))\Big|.\label{eq:upbd1}
\end{multline}
By assumptions $\bold{(H2)}$, $\bold{(H3)}$ and equation \eqref{eq:jacob}, it follows that, up to a positive constant,
\eqref{eq:upbd1} is bounded above by
\[
\sum_{k=1}^n
 {\bf 1}_{\{ \bold X ( D ) = k \}}\frac{\rho(X_1)^{-1}\cdots \rho(X_k)^{-1}}
 {\mathrm{det}\,J[D](
X_1,\ldots,X_{k} ) },
\]
for any $t$ in a neighborhood of zero.
To conclude the proof, we only need to check that the mean of this
r.v. is finite.
We have by definition of the Janossy densities, and since $\mathrm{det}\,J[D](
\bold x ) >0$, for $\PP_D$-a.e. $\bold x \in N_f^D$:
\begin{align}
\mathbb{E}\left[
 {\bf 1}_{\{ \bold X ( D ) = k \}}\frac{\rho(X_1)^{-1} \cdots \rho(X_k)^{-1}}
 {\mathrm{det}\,J[D](
X_1,\ldots,X_{k} ) }\right]
&=\frac{1}{k!} \int_{D^k}\frac{j_{D}^{k}(x_1,\ldots,x_k)}{\mathrm{det}\,J[D](
x_1,\ldots,x_{k} )
}{\bold 1}_{\{j_{D}^{k}(x_1,\ldots,x_k)>0\}}\,\mathrm{d}x_1\ldots\mathrm{d}x_k\nonumber\\
&=\frac{\mathrm{Det}(\bold{Id}-\mathcal{K}_D)}{k!}\,\ell(D^k)<\infty,
\nonumber
\end{align}
where $\ell$ denotes the Lebesgue measure.
\end{Proof}

\begin{remark}
We remark that there is a sign change in (\ref{ibpl}), as compared to the results of \cite{camilier},
which is justified by the corrected formula for (\ref{eq:densitymuphi}). This corrected version is also more in line with the corresponding integration by parts on the Poisson space given in \cite{akr}.
\end{remark}
\noindent
 Next, 
 we extend the integration by parts formula
 by closability to a larger class of functions.
 For $v\in {\cal C}^\infty (D ,\R^d )$, we consider the closability
 of the linear operators
$\nabla_{v}^{\mathrm{N}_{lf}}:\K_D \longrightarrow \mathrm{L}_D^2$ and $\nabla_{v}^{\mathrm{N}_{lf}*}
:\K_D \longrightarrow \mathrm{L}_D^2$. 
 In the
 following, we denote by
 $\overline{A}$ the minimal closed extension
 of a closable linear operator $A$,
 and by
 $\mathrm{Dom}(\overline{A})$ the domain of $\overline{A}$.
\\

\noindent
 In Theorem~\ref{thm:integrationbyparts} below we
 assume, in addition to $\bold{(H1)}$, $\bold{(H2)}$ and $\bold{(H3)}$,
 the following condition.
\\

\noindent
\text{\bf (H4)}: {\it
For any $n\geq 1$, $1\leq i,j\leq n$, and $1\leq h,k\leq d$, we have
\begin{multline*}
\int_{D^n}
 \left|
 \frac{\partial_{x_i^{(h)}}\mathrm{det}\,J[D](x_1,\ldots,x_n)
 \partial_{x_j^{(k)}}\mathrm{det}\,J[D](x_1,\ldots,x_n)}{\mathrm{det}\,J[D](x_1,\ldots,x_n)}
 \right|\\
{\bold 1}_{\left\{ \mathrm{det}\,J[D](x_1,\ldots,x_n) > 0\right\}} \,\mu(\mathrm{d}x_1)\cdots\mu(\mathrm{d}x_n)<\infty.
 \end{multline*}
}
\begin{theorem}
\label{thm:integrationbyparts}
 Assume $\bold{(H1)}-\bold{(H4)}$.

\begin{description}
\item{$(i)$}
 For any $v\in {\cal C}^\infty (D,\R^d )$, the linear operators $\nabla_v^{\mathrm{N}_{lf}}$ and
 $\nabla_v^{\mathrm{N}_{lf}*}$ are well-defined, i.e.
$$
 \nabla_v^{\mathrm{N}_{lf}}(\K_D)\subset\mathrm{L}^2_D \qquad
 \mbox{and}
 \qquad
 \nabla_v^{\mathrm{N}_{lf}*}(\K_D)\subset\mathrm{L}^2_D,
$$
and closable.
\item{$(ii)$}
 For any $v\in {\cal C}^\infty (D,\R^d )$, we have
\[
\mathbb{E}\left[{G}(\bold X^D)\overline{\nabla_v^{\mathrm{N}_{lf}}} {F}(\bold X^D)\right]=
\mathbb{E}\left[{F}(\bold X^D)\overline{\nabla_v^{\mathrm{N}_{lf}*}} {G}(\bold X^D)\right]
\]
 for all
 ${F}\in\mathrm{Dom}\left(\overline{\nabla_v^{\mathrm{N}_{lf}}}\right)$,
 ${G}\in\mathrm{Dom}\left(\overline{\nabla_v^{\mathrm{N}_{lf}*}}\right)$.
\end{description}
\end{theorem}
\noindent Note that under $\bold{(H1)}$,
$\bold{(H2)}$ and $\bold{(H3)}$, condition \text{\bf (H4)} is
satisfied if, for any $n\geq 1$, the function
$$
(x_1,\ldots,x_n)\longmapsto\mathrm{det}\,J[D](x_1,\ldots,x_n),
$$
 is strictly positive on the compact $D^n$.
\begin{Proofy} {\em of Theorem~\ref{thm:integrationbyparts}.}
$(i)$
 Let
 $v\in {\cal C}_c^\infty (D ,\R^d )$ and $F\in{\cal S}_D$. For ease
 of notation, throughout this proof we write $\nabla_v$ in place of
 $\nabla_v^{\mathrm{N}_{lf}}$ and $\nabla_v^{*}$ in place of $\nabla_v^{\mathrm{N}_{lf}*}$.
 We clearly have
$$
 |\nabla_v F(\bold{x})|
 \leq C
$$
 for some constant $C>0$, $\PP_D$-a.e., and
 therefore
 $\nabla_v(\K_D)\subset\mathrm{L}_D^2$.
 The claim $\nabla_v^*(\K_D)\subset\mathrm{L}_D^2$ follows if we
 check that
$\| G (\bold{x}) \nabla_v
U[D](\bold{x})\|_{\mathrm{L}_D^2}<\infty$
 and
$\| G (\bold{x})
B_v^\mu(\bold{x})\|_{\mathrm{L}_D^2}<\infty$
 for any $G \in {\cal S}_D$. The latter relation easily follows noticing that
$$
 |G (\bold{x})
B_v^\mu(\bold{x})|
 \leq C
$$
 for some constant $C>0$, $\PP_D$-a.e..
 Taking
$$
G (\bold{x})
 =
 g_0
 {\bf 1}_{\{ \bold x ( D ) = 0 \}}
 +
 \sum_{k=1}^m
 {\bf 1}_{\{ \bold x ( D ) = k \}}
 g_k( x_1,\ldots,x_k )
$$
 of the form \eqref{F},
 by \eqref{U} we have
\begin{eqnarray*}
 G (\bold{x})
 \nabla_v
 U[D](\bold{x})
 & = &
 -
 \sum_{k=1}^m
 {\bf 1}_{\{ \bold x ( D ) = k \}}
 g_k( x_1,\ldots,x_k )
 \sum_{i=1}^k
 \frac{ \nabla_{x_i}
 \mathrm{det}\,J[D] ( x_1,\ldots,x_k )
}{\mathrm{det}\,J[D] ( x_1,\ldots,x_k ) }
 \cdot
 v(x_i)
,
\end{eqnarray*}
 and for some positive constant $C>0$,
\begin{eqnarray*}
\lefteqn{
 \|
 G
 \nabla_v
 U[D]\|_{\mathrm{L}_D^2}^{2}
}
\\
 & = &
 \sum_{k=1}^m
\frac{1}{k!}  \int_{D^k}
 g_k^2 ( x_1,\ldots,x_k ) {\bold 1}_{\left\{ \mathrm{det}\,J[D](x_1,\ldots,x_k) > 0\right\}}
\\
 & &
 \left(
 \sum_{i=1}^k
 \frac{ \nabla_{x_i}
 \mathrm{det}\,J[D] ( x_1,\ldots,x_k )
}{\mathrm{det}\,J[D] ( x_1,\ldots,x_k )}
 \cdot
 v(x_i)
 \right)^2
 j_{D}^{k}(x_1,\ldots,x_k)
 \mu (\mathrm d x_1) \cdots \mu (\mathrm d x_k)
\\
 & = &
 \mathrm{Det}(\bold{Id}-\mathcal{K}_D)
 \sum_{k=1}^m
 \frac{1}{k!}  \int_{D^k}
 \frac{ g_k^2 ( x_1,\ldots,x_k ) }{\mathrm{det}\,J[D] ( x_1,\ldots,x_k )}
\\
 & &
 {\bold 1}_{\left\{ \mathrm{det}\,J[D](x_1,\ldots,x_k) > 0\right\}}
 \left(
 \sum_{i=1}^k
 \nabla_{x_i}
 \mathrm{det}\,J[D] ( x_1,\ldots,x_k )
 \cdot
 v(x_i)
 \right)^2
 \mu (\mathrm dx_1) \cdots \mu (\mathrm dx_k)
\\
 & \leq &
C\, \mathrm{Det}(\bold{Id}-\mathcal{K}_D)
 \sum_{k=1}^m
 \frac{1}{k!}
  \sum_{1 \leq i,j \leq k}
 \int_{D^k}
 {\bold 1}_{\left\{ \mathrm{det}\,J[D](x_1,\ldots,x_k) > 0\right\}}
\\
 & &
 \frac{
 \nabla_{x_i}
 \mathrm{det}\,J[D] ( x_1,\ldots,x_k )
 \cdot
 v(x_i)\,
\nabla_{x_j}
 \mathrm{det}\,J[D] ( x_1,\ldots,x_k )
 \cdot
 v(x_j)}{\mathrm{det}\,J[D] ( x_1,\ldots,x_k )}\,
 \mu (\mathrm dx_1) \cdots \mu (\mathrm dx_k)
\\
 & < & \infty,
\end{eqnarray*}
 where the latter integral is finite by
\text{\bf (H4)}.
\\

\noindent
 To conclude, we only need to
 show that $\nabla_v$ is closable (the closability of
 $\nabla_v^*$ can be proved similarly). Let $(F_n)_{n\geq
1}$ be a sequence in $\K_D$ converging to $0$ in $\mathrm{L}^{2}_D$
and such that $\nabla_v F_n$ converges to $V$ in $\mathrm{L}^{2}_D$
 as $n$ goes to infinity.
 We need to show that $V=0$ $\PP_D$-a.e..
 We have
\begin{align}
\left|\langle G,V\rangle_{\mathrm L^2_D}\right|&=\lim_{n\to\infty}|\mathbb{E}[G(\bold{X}^D)\nabla_v
F_n(\bold{X}^D)]|
=\lim_{n\to\infty}|\mathbb{E}[
F_n(\bold{X}^D)
\nabla_v^* G (\bold{X}^D
]|\label{eq:applintpart1}\\
&\leq\|\nabla_v^* G
\|_{\mathrm{L}_D^2}\lim_{n\to\infty}\|F_n\|_{\mathrm{L}_D^2}=0,\qquad
 G\in{\cal S}_D.
\nonumber
\end{align}
\noindent
 Here, the second inequality in \eqref{eq:applintpart1}
 follows by the integration by parts formula
 \eqref{ibpl}. 
 The conclusion follows noticing that $\langle G,V\rangle_{\mathrm L^2_D}=0$ for all
 $G\in{\cal S}_D$
 implies $V=0$ $\PP_D$-a.e. due to the density
 of $\mathcal S_D$ in $\mathrm L^2_D$.
\\

\noindent
$(ii)$
 By $(i)$ 
 both operators $\nabla_v$
and $\nabla_v^*$ are closable. Take $
F\in\mathrm{Dom}(\overline{\nabla}_v)$, $
G\in\mathrm{Dom}(\overline{\nabla_v^*})$ and let
$(F_n)_{n\geq 1}$, $(G_n)_{n\geq 1}$ be
sequences in $\K_D$ such that $F_n$ converges
 to $ F$, $G_n$ converges to $ G$,
 $\nabla_v F_n$ converges to $\overline{\nabla_v} F$
 and $\nabla_v^* G_n$ converges to $\overline{\nabla_v^*}
{G}$ in $\mathrm{L}_D^2$ as $n$ goes to infinity.
 By Lemma~\ref{le:integrationbyparts}
the integration by parts formula applies to r.v.'s in
$\mathcal{S}_D$, therefore we have $\mathbb{E}[G_n
(\bold{X}^D)\nabla_v F_n(\bold{X}^D)]=\mathbb{E}[F_n(\bold{X}^D)
\nabla_v^* G_n(\bold{X}^D)]$ for all $n\geq 1$. The claim follows
if we prove
 that
\[
\lim_{n\to\infty}\mathbb{E}[G_n(\bold{X}^D)\nabla_v F_n(\bold
{X}^D)]=\mathbb{E}[ G(\bold X^D)\overline{\nabla_v} F(\bold X^D)]
\]
and
\[
\lim_{n\to\infty}\mathbb{E}[F_n(\bold{X}^D)\nabla_v^*
G_n(\bold{X}^D)]=\mathbb{E}[ F(\bold X^D)\overline{\nabla_v^*}
G(\bold X^D)].
\]
We only show the first limit above; the second limit being proved
similarly. We have
\begin{align}
&|\mathbb{E}[G_n(\bold{X}^D)\nabla_v F_n(\bold
{X}^D)]-\mathbb{E}[
G(\bold X^D)\overline{\nabla_v} F(\bold X^D)]|\nonumber\\
&=|\mathbb{E}[G_n(\bold{X}^D)\nabla_v F_n(\bold
{X}^D)]-\mathbb{E}[G_n(\bold{X}^D)\overline{\nabla_v}
F(\bold X^D)]+\mathbb{E}[G_n(\bold{
X}^D)\overline{\nabla_v} F(\bold X^D)]\nonumber\\
&\qquad\qquad\qquad\qquad\qquad\qquad\qquad\qquad\qquad\qquad\qquad-\mathbb{E}[ G(\bold X^D)\overline{\nabla_v} F(\bold X^D)]|\nonumber\\
&\leq|\mathbb{E}[G_n(\bold{X}^D)(\nabla_v F_n(\bold
{X}^D)-\overline{\nabla_v} F(\bold X^D))]|+|\mathbb{E}[(G_n(\bold
{X}^D)-{G}(\bold X^D))\overline{\nabla_v}
 F(\bold X^D)]|\nonumber\\
&\leq\|G_n\|_{\mathrm{L}_D^2}\|\nabla_v
F_n-\overline{\nabla_v}{F}\|_{\mathrm{L}_D^2}
+\|G_n-
G\|_{\mathrm{L}_D^2}\|\overline{\nabla_v}
F\|_{\mathrm{L}_D^2} , \nonumber
\end{align}
which tends to $0$ as $n$ goes to infnity.
\end{Proofy}
\bigskip
\noindent
\section{Dirichlet forms corresponding to determinantal processes on $D\subset S$} 
\label{sec:df}
 In this section we construct the Dirichlet form
 associated to a determinantal process,
 cf. Theorem~\ref{thm:dirichlet} below.
 We start by recalling some definitions related to bilinear forms (see
\cite{rockner} for details). Let $H$ be a Hilbert space with inner
product $\langle \cdot,\cdot \rangle$ and
$\mathcal{A}:\mathrm{Dom}(\mathcal{A})\times\mathrm{Dom}(\mathcal{A})
  \longrightarrow \R$ a bilinear form defined on a dense subspace
$\mathrm{Dom}(\mathcal{A})$ of $H$, the domain of $\mathcal A$.
The form $\mathcal{A}$ is said to be symmetric if
$\mathcal{A}(F,G)=\mathcal{A}(G,F)$, for any
$F,G\in\mathrm{Dom}(\mathcal{A})$, and 
 non-negative definite if $\mathcal{A}(F,F)\geq 0$, for any
$F\in\mathrm{Dom}(\mathcal{A})$. Let $\mathcal{A}$ be symmetric
and non-negative definite, $\mathcal A$ is said to be closed if
$\mathrm{Dom}(\mathcal A)$ equipped with the norm
\[
\|F\|_{\mathcal{A}}:=\sqrt{\mathcal{A}(F,F)+ \langle F,F \rangle },\quad\text{$F\in\mathrm{Dom}(\mathcal{A})$},
\]
is a Hilbert space. A symmetric and non-negative definite bilinear
form $\mathcal{A}$ is said to be closable if, for any sequence
$(F_n)_{n\geq 1}\subset\mathrm{Dom}(\mathcal A)$ such that $F_n$
goes to $0$ in $H$ and $(F_n)_{n\geq 1}$ is Cauchy {\em w.r.t.}
$\|\cdot\|_\mathcal{A}$ it holds that $\mathcal{A}(F_n,F_n)$
converges to $0$ in $\R$ as
 $n$ goes to infinity.
 Let $\mathcal{A}$ be
closable and denote by $\mathrm{Dom}(\overline{\mathcal{A}})$ the
completion of $\mathrm{Dom}(\mathcal{A})$ {\em w.r.t.} the norm
$\|\cdot\|_\mathcal{A}$. It turns out that $\mathcal A$ is
uniquely extended to $\mathrm{Dom}(\overline{\mathcal{A}})$ by the
closed, symmetric and non-negative definite bilinear form
\[
\overline{\mathcal
A}(F,G)=\lim_{n\to\infty}\mathcal{A}(F_n,G_n),\quad\text{$(F,G)\in\mathrm{Dom}(\overline{\mathcal{A}})\times
\mathrm{Dom}(\overline{\mathcal{A}})$,}
\]
 where $\{(F_n,G_n)\}_{n\geq 1}$ is any sequence in
$\mathrm{Dom}(\mathcal{A})\times\mathrm{Dom}(\mathcal{A})$
such that $(F_n,G_n)$ converges to $(F,G)\in
\mathrm{Dom}(\overline{\mathcal{A}})\times
\mathrm{Dom}(\overline{\mathcal{A}})$ {\em w.r.t.} the norm
$\|\cdot\|_{\overline{\mathcal
A}}+\|\cdot\|_{\overline{\mathcal{A}}}$. 
Suppose $H=\mathrm{L}^2(B,\beta)$ where $(B,\mathcal F,\beta)$
is a $\sigma$-finite measure space. 
A symmetric, non-negative definite and closed bilinear 
form $\mathcal A$ is said to be a Dirichlet form if
\[
\mathcal{A}(\min\{F^+,1\},\min\{F^+,1\})\leq\mathcal{A}(F,F),
\qquad\text{$F\in \mathrm{Dom}(\mathcal A)$,}
\]
where $F^+$ denotes the positive part of $F$. Suppose that $B$ is
a Hausdorff topological space, $\mathcal F=\mathcal B(B)$ the corresponding Borel $\sigma$-algebra, and let $\mathcal A$ be a
Dirichlet form.
Throughout this paper, we shall use the following notions.
The generator $\mathcal H^{\mathrm{gen}}$ of the Dirichlet form is the unique symmetric non-negative definite operator on $H$ such that
\begin{equation*}
	\mathcal A(F,G)=\langle\left(-\mathcal H^{\mathrm{gen}}F\right),G\rangle,
\end{equation*}
for $F,G\in\mathrm{Dom}(-\mathcal H^{\mathrm{gen}})\subset\mathrm{Dom}(\mathcal A)$.
The symmetric semi-group of
$\mathcal A$
is the linear operator on $H$ defined by $T_t
F:=\mathrm{exp}(-t\mathcal H^{\mathrm{gen}}) F$,
$t>0$. 
An $\mathcal{A}$-nest is an increasing sequence
$(C_n)_{n\geq 1}$ of closed subsets of $B$ such that
\[
\text{$\bigcup_{n\geq 1}$\{$F\in\mathrm{Dom}(\mathcal A)$: $F=0$
$\beta$-a.e. on $B\setminus C_n$\}}
\]
is dense in $\mathrm{Dom}(\mathcal A)$ {\em w.r.t.} the norm
$\|\cdot\|_{\mathcal{A}}$. Throughout this paper, we say that a subset $B'\subset B$ is
$\mathcal{A}$-exceptional if there exists an $\mathcal A$-nest
$(C_n)_{n\geq 1}$ with $B' \subset B \setminus \bigcup_{n\geq 1}C_n$.
We say that a property holds
$\mathcal{A}$-almost everywhere ($\mathcal{A}$-a.e.) if it holds
up to an $\mathcal A$-exceptional set. We say that a function
$f:B\to\R$ is $\mathcal{A}$-almost continuous
($\mathcal{A}$-a.c.) if there exists an $\mathcal A$-nest
$(C_n)_{n\geq 1}$ such that the restriction $f_{|C_n}$ of $f$ to
$C_n$ is continuous for each $n\geq 1$.
We say that a Dirichlet
form $\mathcal{A}$ 
is quasi-regular if:

\noindent $(i)$ There exists an $\mathcal{A}$-nest $(C_n)_{n\geq
1}$ consisting of compact sets.
\\
$(ii)$ There exists a $\|\cdot\|_{\mathcal{A}}$-dense subset of
$\mathrm{Dom}(\mathcal{A})$ whose elements have $\mathcal{A}$-a.c.
$\beta$-versions.
\\
$(iii)$ There exist $F_k\in\mathrm{Dom}(\mathcal{A})$, $k\geq 1$,
having $\mathcal{A}$-a.c. $\beta$-versions $\tilde{F}_k$, $k\geq
1$, such that $(\tilde{F}_k)_{k\geq 1}$ is a separating set for
$B\setminus N$ (i.e. for any $x,y\in B\setminus N$, $x\neq y$,
there exists ${k^*}\ge1$ such that
$\tilde{F}_{k^*}(x)\neq\tilde{F}_{k^*}(y)$),
where $N$ is a subset of $B$ which is $\mathcal
A$-exceptional.

After these general considerations, we move to the situation at hand.
Assume \text{$\bold{(H1)}$},  \text{$\bold{(H2)}$} and \text{$\bold{(H3)}$}.
In particular, we recall that $D$ is a fixed compact set of $S$ which is in turn included in $\R^d$.
We denote by
$\mathrm{\ddot{N}}_{f}^{D}$ the set of $\mathbb N$-valued
Radon measures on $D$. We equip $\mathrm{\ddot{N}}_{f}^{D}$ with the vague topology and the corresponding Borel $\sigma$-algebra $\mathcal{B}(\mathrm{\ddot{N}}_{f}^D)$.
Note that $\mathrm{N}_f^D$ is contained in $\mathrm{\ddot{N}}_{f}^D$.
In the following, we consider the subspace $\tilde{\cal S}_D$ of
${\cal S}_D$ made of cylindrical functions.
\begin{definition}
\label{def:testfunctions2}
A function  $F:\mathrm N_f^D\to\R$ is said to be in $\tilde{\cal S}_D$ if it is of the
form
\begin{equation*}
\label{eq:deftildeS}
F(\bold{x})=f\left(\sum_{k=1}^{\bold x (D)}\varphi_1(x_k),\ldots,\sum_{k=1}^{\bold x (D)}\varphi_m(x_k)\right)
{\bf 1}_{\{ \bold x (D) \leq n \}},\qquad\bold x=\{x_1,\dots,x_{\mathbf x(D)}\}\in\mathrm N_f^D,
\end{equation*}
for some integers $m,n \geq 1$, $\varphi_1,\ldots,\varphi_m \in
{\cal C}^\infty(D)$, $f\in {\cal C}_b^\infty(\R^m)$.
\end{definition}
Note that, as already mentioned, $\tilde{\cal S}_D$ is dense in $\mathrm{L}_D^2$.
We consider the bilinear map $\mathcal{E}_D$ defined on
$\tilde{\mathcal{S}}_D\times\tilde{\mathcal{S}}_D$ by
\begin{equation*}
\mathcal{E}_D(F,G)
:=\mathbb{E}\left[\sum_{i=1}^{\bold X (D)}\nabla_{X_i}^{\mathrm{N}_{lf}}F(\bold{X}^D)\cdot\nabla_{X_i}^{\mathrm{N}_{lf}}G(\bold{X}^D)\right].
\end{equation*}
For $F\in {\cal S}_D$ of the form \eqref{F}, i.e.
$$
 F (\bold{x} )
 =
 f_0
 {\bf 1}_{\{ \bold x ( D ) = 0 \}}
 +
 \sum_{k=1}^n
 {\bf 1}_{\{ \bold x ( D ) = k \}}
 f_k( x_1,\ldots,x_k ),\qquad\bold x=\{x_1,\dots,x_{\mathbf x(D)}\}\in\mathrm N_f^D
$$
 we also define the Laplace operator $\mathcal{H}_D$ by
\begin{eqnarray*}
\lefteqn{
 \mathcal{H}_DF(\bold{x} ) =
 \sum_{k=1}^n
 {\bf 1}_{\{ \bold x ( D ) = k \}}
}
\\
 & &
 \! \! \! \! \! \! \! \! \! \! \! \!
 \sum_{i=1}^k
 \left(-
 \beta^\mu(x_i)\cdot\nabla_{x_i}  f_k( x_1,\ldots,x_k )
 -
 \Delta_{x_i} f_k( x_1,\ldots,x_k )
 +
 U_{i,k}( x_1,\ldots,x_k )
 \cdot \nabla_{x_i} f_k ( x_1,\ldots,x_k )
 \right),
\end{eqnarray*}
 where $\Delta=-\mathrm{div}\nabla$ denotes the Laplacian.
The next theorem provides the Dirichlet form
associated to a determinantal process.
\begin{theorem}\label{thm:dirichlet}
 Under conditions $\bold{(H1)}-\bold{(H4)}$
 we have:
\\
\noindent $(i)$ The Laplace operator $\mathcal{H}_D:\tilde{\cal S}_D
\longrightarrow \mathrm{L}_D^2$ is linear, symmetric, non-negative
definite and well-defined, i.e. $\mathcal{H}_D(\tilde{\cal
S}_D)\subset\mathrm{L}_D^2$. In particular the operator square
root $\mathcal{H}_D^{1/2}$ of $\mathcal{H}_D$ exists.
\\
\noindent $(ii)$ The bilinear form $\mathcal{E}_D:\tilde{\cal
S}_D\times\tilde{\cal S}_D \longrightarrow \R$ is symmetric, non-negative
definite and well-defined, i.e. $\mathcal{E}_D(\tilde{\cal
S}_D\times\tilde{\cal
S}_D)\subset\R$.\\
\noindent $(iii)$ $\mathcal{H}_D^{1/2}$ and $\mathcal{E}_D$ are
closable and the following relation holds:
\begin{equation}\label{eq:Dirichlet}
\overline{\mathcal{E}_D}({F},{G})=\mathbb{E}\left[\overline{\mathcal{H}_D^{1/2}}\,{F}(\mathbf X^D)\,
\overline{\mathcal{H}_D^{1/2}}\,{G}(\mathbf X^D)\right], \qquad\text{$\forall$
${F},{G}\in\mathrm{Dom}(\overline{\mathcal{H}_D^{1/2}})$.}
\end{equation}
\noindent $(iv)$ The bilinear form
$(\overline{\mathcal{E}_D},\mathrm{Dom}(\overline{\mathcal{H}_D^{1/2}}))$
is a Dirichlet form.
\end{theorem}
The proof of the theorem is based on the following lemma.
\begin{lemma}\label{le:dirichlet}
 Under conditions $\bold{(H1)}-\bold{(H4)}$,
 for any $F,G\in \tilde{\cal S}_D$, we have
\begin{align}
\mathbb{E}\left[\sum_{i=1}^{\bold X (D)}\nabla_{X_i}^{\mathrm{N}_{lf}}F(\bold{X}^D)\cdot\nabla_{X_i}^{\mathrm{N}_{lf}}G(\bold{X}^D)\right]
&=\mathbb{E}[G(\bold{X}^{D})\mathcal{H}_DF(\bold{X}^{D})]\label{eq:Dir1}\\
&=\mathbb{E}[\mathcal{H}_D^{1/2}F(\bold{X}^{D})\mathcal{H}_D^{1/2}G(\bold{X}^{D})].\label{eq:Dir2}
\end{align}
\end{lemma}
\begin{Proofy} {\em of Theorem~\ref{thm:dirichlet}.}
 
\noindent$\it{(i)}$ By Relation~\eqref{eq:Dir1} in Lemma
\ref{le:dirichlet} we easily deduce that, for any
$F,G\in \tilde{\cal S}_D$ we have
\[
\mathbb{E}[G(\bold{X}^D)\mathcal{H}_DF(\bold{X}^D)]=\mathbb{E}[F(\bold{X}^D)\mathcal{H}_DG(\bold{X}^D)]
\quad\text{and}\quad
\mathbb{E}[F(\bold{X}^D)\mathcal{H}_DF(\bold{X}^D)]\geq 0.
\]
Therefore, $\mathcal{H}_D$ is symmetric and non-negative definite.
It remains to check that, under the foregoing assumptions,
$\mathcal{H}_D$ is well-defined. Let $F\in\tilde{\mathcal{S}}_D$ be of the form
\begin{align}
F(\bold{x})&=\sum_{k=1}^{n}{\bf 1}_{\{ \bold x (D) = k \}}f\left(\sum_{i=1}^{k}\varphi_1(x_i),\ldots,\sum_{i=1}^{k}\varphi_m(x_i)\right)\nonumber\\
&=\sum_{k=1}^{n}{\bf 1}_{\{ \bold x (D) = k \}}f_k(x_1,\ldots,x_k)\nonumber
\end{align}
for some integers $m,n \geq 1$, $\varphi_1,\ldots,\varphi_m \in
{\cal C}^\infty(D)$, $f\in {\cal C}_b^\infty(\R^m)$.
For the well-definiteness of $\mathcal{H}_D$ we only need to check that, for $\mathbf x=\{x_1,\dots,x_{\mathbf x(D)}\}\in\mathrm N_f^D$,
$$F_1(\bold x):={\bold 1}_{\{\mathbf x(D)=k\}}\left(\beta^\mu(x_i)\right)^{(j)},\text{ and }F_2(\bold x):={\bold 1}_{\{\mathbf x(D)=k\}}\left(U_{i,k}(x_1,\dots,x_k)\right)^{(j)},$$
are in $\mathrm L^2_D$ for all $k\in\{1,\dots,n\}$, $i\in\{1,\dots,k\}$ and $j\in\{1,\dots,d\}$.
One may easily check that $F_1\in\mathrm L^2_D$ due to \text{$\bold{(H2)}$} and $F_2\in\mathrm L^2_D$ due to \text{$\bold{(H4)}$}.
\\

\noindent$\it{(ii)}$
 The symmetry and non-negative definiteness of
 $\mathcal{E}_D$ 
 follow from Lemma~\ref{le:dirichlet}.
 It remains to check that, under the
foregoing assumptions, $\mathcal{E}_D$ is well-defined. By Step
$(i)$, for any $F\in \tilde{\cal S}_D$, we have
$\mathcal{H}_DF\in\mathrm{L}_D^2$. We conclude the proof by
noting that, by Lemma~\ref{le:dirichlet}, for any
$F,G\in \tilde{\cal S}_D$ and some
positive constant $c>0$, we have
\[
|\mathcal{E}_D(F,G)|=|\mathbb{E}[G(\bold{X}^D)\mathcal{H}_DF(\bold{X}^D)]|
\leq c\|\mathcal{H}_DF(\bold{X}^D)\|_{\mathrm{L}_D^2}<\infty.
\]

\noindent$\it{(iii)}$
 We first show that $\mathcal{E}_D$ is closable. By Lemma 3.4 page 29 in \cite{rockner}, we have to check that 
if $(F_n)_{n\geq 1}\subset
\tilde{\cal S}_D$ is such that $F_n$ converges to $0$
in $\mathrm{L}_D^2$, then
$\mathcal{E}_D(G,F_n)$ converges to $0$, for
any $G\in \tilde{\cal S}_D$. This easily follows by
Lemma~\ref{le:dirichlet}, the Cauchy-Schwarz inequality and the
fact that $\mathcal{H}_DG$ is square integrable (see Step $(i)$). The closability of $\mathcal{H}_D^{1/2}$
follows by the closability of $\mathcal E_D$,
relation~\eqref{eq:Dir2} in Lemma~\ref{le:dirichlet} and Remark
3.2 $(i)$ page~29 in \cite{rockner}. Finally, we prove relation~\eqref{eq:Dirichlet}. Take
$ F, G\in\mathrm{Dom}(\overline{\mathcal{H}_D^{1/2}})$ and
let $(F_n)_{n\geq 1}$, $(G_n)_{n\geq 1}$
be sequences in $\tilde{\mathcal{S}}_D$ such that
$F_n$ converges to $ F$, $G_n$
 converges to $ G$,
$\mathcal{H}_D^{1/2}
F_n$
 converges to $\overline{\mathcal{H}_D^{1/2}}\, F$,
 and
$\mathcal{H}_D^{1/2}G_n$
 converges to $\overline{\mathcal{H}_D^{1/2}}\,{G}$
in $\mathrm{L}_D^2$ as $n$ goes to infinity.
 By Lemma~\ref{le:dirichlet} we have
\[
\mathcal{E}_D(F_n,G_n
)=\mathbb{E}[\mathcal{H}_D^{1/2}
F_n(\bold{X}^D)\mathcal{H}_D^{1/2} G_n(\bold{X}^D)],\quad\text{for
all $n\geq 1$.}
\]
The claim follows if we prove that 
\[
\lim_{n\to\infty}\mathbb{E}[\mathcal{H}_D^{1/2}F_n(\bold
{X}^D)\mathcal{H}_D^{1/2}G_n(\bold
{X}^D)]=\mathbb{E}[\overline{\mathcal{H}_D^{1/2}}\,
F(\bold
{X}^D)\overline{\mathcal{H}_D^{1/2}}\, G(\bold
{X}^D)].
\]
Indeed we have:
\begin{align}
&|\mathbb{E}[\mathcal{H}_D^{1/2}G_n(\bold{X}^D)\mathcal{H}_D^{1/2}F_n(\bold
{X}^D)]-\mathbb{E}[\overline{\mathcal{H}_D^{1/2}}
G(\bold
{X}^D)\,\overline{\mathcal{H}_D^{1/2}}\, F(\bold
{X}^D)]|\nonumber\\
&=|\mathbb{E}[\mathcal{H}_D^{1/2}G_n(\bold{X}^D)\mathcal{H}_D^{1/2}F_n(\bold
{X}^D)]-\mathbb{E}[\mathcal{H}_D^{1/2}G_n(\bold{X}^D)\,\overline{\mathcal{H}_D^{1/2}}\,
F(\bold
{X}^D)]\nonumber\\
&\,\,\,\,\,\,\,\,\,\,\,\,\,\,\,\,\,\,\,\,\,\,\,\,\,\,\,\,\,\,\,\,\,\,\,\,
+\mathbb{E}[\mathcal{H}_D^{1/2}G_n(\bold{
X}^D)\,\overline{\mathcal{H}_D^{1/2}}\,
F(\bold
{X}^D)]-\mathbb{E}[\overline{\mathcal{H}_D^{1/2}} G(\bold
{X}^D)\,
\overline{\mathcal{H}_D^{1/2}}\, F(\bold
{X}^D)]|\nonumber\\
&\leq|\mathbb{E}[\mathcal{H}_D^{1/2}G_n(\bold{X}^D)(\mathcal{H}_D^{1/2}F_n(\bold
{X}^D)-\overline{\mathcal{H}_D^{1/2}}\,
F(\bold
{X}^D))]|\nonumber\\
&\,\,\,\,\,\,\,\,\,\,\,\,\,\,\,\,\,\,\,\,\,\,\,\,\,\,\,\,\,\,\,\,\,\,\,\,+|\mathbb{E}[(\mathcal{H}_D^{1/2}G_n(\bold
{X}^D)-\overline{\mathcal{H}_D^{1/2}}{G}(\bold
{X}^D))\overline{\mathcal{H}_D^{1/2}}\,
 F(\bold
{X}^D)]|\nonumber\\
&\leq\|\mathcal{H}_D^{1/2}G_n\|_{\mathrm{L}_D^2}\|\mathcal{H}_D^{1/2}
F_n-\overline{\mathcal{H}_D^{1/2}}\,{F}\|_{\mathrm{L}_D^2}\nonumber\\
&\,\,\,\,\,\,\,\,\,\,\,\,\,\,\,\,\,\,\,\,\,\,\,\,\,\,\,\,\,\,\,\,\,\,\,\,
+\|\mathcal{H}_D^{1/2}G_n-\overline{\mathcal{H}_D^{1/2}}
G\|_{\mathrm{L}_D^2}\|\overline{\mathcal{H}_D^{1/2}}\,
F\|_{\mathrm{L}_D^2}\to 0,\quad\text{as $n\to\infty$.}\nonumber
\end{align}

\noindent$\it{(iv)}$
The bilinear form
$(\overline{\mathcal{E}_D},\mathrm{Dom}(\overline{\mathcal{H}_D^{1/2}}))$
defined by \eqref{eq:Dirichlet} is clearly symmetric, non-negative
definite, and closed. 
We conclude the
proof by applying Proposition 4.10 page~35 in \cite{rockner}. First, 
note that $\tilde{\mathcal{S}}_D$ is dense in
$\mathrm{Dom}(\overline{\mathcal{H}_D^{1/2}})$ ({\em w.r.t.} the
norm $\|\cdot\|_{\overline{\mathcal E}_D}$). By Exercise 2.7 page~47
in \cite{rockner}, for any $\varepsilon>0$ there exists an
infinitely differentiable function $\varphi_\varepsilon:\R
 \longrightarrow
[-\varepsilon,1+\varepsilon]$ (which shall not be confused with the functions
$\varphi_1,\ldots,\varphi_m$ involved in Definition~\ref{def:testfunctions2}) such that $\varphi_\varepsilon(t)=t$
for any $t\in [0,1]$,
$0\leq\varphi_\varepsilon(t)-\varphi_\varepsilon(s)\leq t-s$ for
all $t,s\in\R$, $t\geq s$, $\varphi_\varepsilon(t)=1+\varepsilon$
for $t\in [1+2\varepsilon,\infty)$ and
$\varphi_\varepsilon(t)=-\varepsilon$ for $t\in
(-\infty,-2\varepsilon]$. Note that
$|\varphi_\varepsilon^\prime(t)|^2\leq 1$ for any $\varepsilon>0$,
$t\in\R$ and $\varphi_\varepsilon$ is in ${\cal C}_b^\infty$, for
any $\varepsilon>0$. Consider the function
\[
F(\bold{x})=f\left(\sum_{k=1}^{\bold x (D)}\varphi_1(x_k),\ldots,\sum_{k=1}^{\bold x (D)}\varphi_m(x_k)\right)
{\bf 1}_{\{ \bold x (D) \leq n \}},\qquad\mathbf x=\{x_1,\dots,x_{\mathbf x(D)}\}\in\mathrm N_f^D,
\]
for some integers $m,n \geq 1$, $\varphi_1,\ldots,\varphi_m \in
{\cal C}^\infty(D)$, $f\in {\cal C}_b^\infty(\R^m)$. Note that
$\varphi_{\varepsilon}\circ
F\in\tilde{\mathcal{S}}_D$. Indeed we have 
\begin{align}
\varphi_{\varepsilon}\circ F(\bold{x})&=
\varphi_{\varepsilon}\left(f\left(\sum_{k=1}^{\bold x (D)}\varphi_1(x_k),\ldots,\sum_{k=1}^{\bold x (D)}\varphi_m(x_k)\right)
{\bf 1}_{\{ \bold x (D) \leq n \}}\right)\nonumber\\
&=\varphi_\varepsilon\left(f\left(\sum_{k=1}^{\bold x (D)}\varphi_1(x_k),\ldots,\sum_{k=1}^{\bold x (D)}\varphi_m(x_k)\right)\right)
{\bf 1}_{\{ \bold x (D) \leq n \}},\nonumber
\end{align}
because $\varphi_\varepsilon(0)=0$. Next, we have
\begin{align}
&\overline{\mathcal{E}_D}(\varphi_\varepsilon\circ
F,\varphi_\varepsilon\circ
F)=\mathbb{E}\left[\sum_{i=1}^{\bold X (D)}\nabla_{X_i}^{\mathrm{N}_{lf}}\varphi_\varepsilon\circ
F(\bold{X}^D)\cdot\nabla_{X_i}^{\mathrm{N}_{lf}}\varphi_\varepsilon\circ
F(\bold{X}^D)\right]\nonumber\\
&=\mathbb{E}\left[\sum_{i=1}^{\bold X (D)}\|\nabla_{X_i}^{\mathrm{N}_{lf}}\varphi_\varepsilon\circ
F(\bold{X}^D)\|^2\right]\nonumber\\
&=\mathbb{E}\left[\sum_{k=1}^{n}{\bf 1}_{\{ \bold X (D) = k
\}}\sum_{i=1}^{k}\|\nabla_{X_i}^{\mathrm{N}_{lf}}\varphi_\varepsilon\circ
F(\bold{X}^D)\|^2\right]\nonumber\\
&=\mathbb{E}\Biggl[\sum_{k=1}^{n}{\bf 1}_{\{ \bold X (D) = k
\}}\sum_{i=1}^{k}\Big\|\sum_{j=1}^{m}
\varphi_\varepsilon^{\prime}\circ
f\left(\sum_{l=1}^{k}\varphi_1(X_l),\ldots,\sum_{l=1}^{k}\varphi_m(X_l)\right)
\nonumber\\
&\,\,\,\,\,\,\,\,\,\,\,\,\,\,\,\,\,\,\,\,\,\,\,\,\,\,\,\,\,\,\,\,\,\,\,\,\,\,\,\,
\,\,\,\,\,\,\,\,\,\,\,\,\,\,\,\,\,\,\,\,\,\,\,\,\,\,\,\,\,\,\,\,\,\,\,\,\,\,\,\,
\times
\partial_{j}
f\left(\sum_{l=1}^{k}\varphi_1(X_l),\ldots,\sum_{l=1}^{k}\varphi_m(X_l)\right)
\nabla\varphi_j(X_i)\Big\|^2\Biggr]\nonumber\\
&\leq\mathbb{E}\Biggl[\sum_{k=1}^{n}{\bf 1}_{\{ \bold X (D) = k
\}}\sum_{i=1}^{k}\Big\|\sum_{j=1}^{m}\partial_{j}
f\left(\sum_{l=1}^{k}\varphi_1(X_l),\ldots,\sum_{l=1}^{k}\varphi_m(X_l)\right)
\nabla\varphi_j(X_i)\Big\|^2\Biggr]\label{eq:Dir3}\\
&=\overline{\mathcal{E}_D}(F,F),\nonumber
\end{align}
 where in \eqref{eq:Dir3} we used the fact that
 $| \varphi_\varepsilon^\prime(t) |^2 \leq 1$, $t\in\R$. By this
 inequality we have, for any $F\in \tilde{\cal S}_D$,
\[
\liminf_{\varepsilon\to
0}\overline{\mathcal{E}_D}(F\pm\varphi_\varepsilon\circ
F,F\mp\varphi_\varepsilon\circ
F)\geq 0
\]
and the proof is completed (since, as required by Proposition 4.10
page~35 in \cite{rockner}, we checked condition (4.6) page~34 in
\cite{rockner}). Indeed, for any $\varepsilon>0$, by the above
inequality and Lemma~\ref{le:dirichlet}, we
have
\begin{align}
&\overline{\mathcal{E}_D}(F+\varphi_\varepsilon\circ
F,F-\varphi_\varepsilon\circ
F)\nonumber\\
&=\overline{\mathcal{E}_D}(F-\varphi_\varepsilon\circ
F,F+\varphi_\varepsilon\circ
F)\nonumber\\
&=\mathbb{E}[(F(\bold{X}^D)-\varphi_\varepsilon\circ
F(\bold{X}^D))\mathcal{H}_D(F(\bold{X}^D)+\varphi_\varepsilon\circ
F(\bold{X}^D))]\nonumber\\
&=\mathbb{E}[F(\bold{X}^D)\mathcal{H}_DF(\bold{X}^D)+
F(\bold{X}^D)\mathcal{H}_D\varphi_\varepsilon\circ
F(\bold{X}^D)\nonumber\\
&\,\,\,\,\,\,\,\,\,\,\,\,-\varphi_\varepsilon\circ
F(\bold{X}^D)\mathcal{H}_DF(\bold{X}^D)-\varphi_\varepsilon\circ
F(\bold{X}^D)\mathcal{H}_D\varphi_\varepsilon\circ
F(\bold{X}^D)]
\nonumber\\
&\geq\mathbb{E}[F(\bold{X}^D)\mathcal{H}_D\varphi_\varepsilon\circ
F(\bold{X}^D)-\varphi_\varepsilon\circ
F(\bold{X}^D)\mathcal{H}_DF(\bold{X}^D)]
\nonumber\\
&=0.\nonumber
\end{align}
\end{Proofy}

\begin{Proofy} {\em of Lemma~\ref{le:dirichlet}.}
 Let
 $F,G\in \tilde{\cal S}_D$ be, respectively, defined for $\mathbf x=\{x_1,\dots,x_{\mathbf x(D)}\}\in\mathrm N_f^D$ as
\[
 F(\bold{x})=f\left(\sum_{k=1}^{\bold x (D)}\varphi_1(x_k),\ldots,\sum_{k=1}^{\bold x (D)}\varphi_{m_1}(x_k)\right)
 {\bf 1}_{\{ \bold x (D) \leq n_1 \} },
\]
\[
G(\bold{x})=g\left(\sum_{k=1}^{\bold x (D)}\gamma_1(x_k),\ldots,\sum_{k=1}^{\bold x (D)}\gamma_{m_2}(x_k)\right)
 {\bf 1}_{\{ \bold x (D) \leq n_2 \} },
\]
for some integers $m_1,m_2,n_1,n_2 \geq 1$,
 $\varphi_1,\ldots,\varphi_{m_1},\gamma_1,\ldots,\gamma_{m_2}
 \in
 {\cal C}^\infty(D)$, $f\in {\cal C}_b^\infty(\R^{m_1})$, $g\in {\cal C}_b^\infty(\R^{m_2})$.
 Define
\[
F_i(\bold{x})=\partial_i
f\left(\sum_{k=1}^{\bold x (D)}\varphi_1(x_k),\ldots,\sum_{k=1}^{\bold x (D)}\varphi_{m_1}(x_k)\right)
 {\bf 1}_{\{ \bold x (D) \leq n_1 \} },
\]
and
\[
v_i(x)=\nabla\varphi_i (x),\quad\text{$x\in D$.}
\]
By direct computation we find
\begin{eqnarray*}
\lefteqn{
 \mathcal{H}_DF(\bold{x} )
 =
 -{\bf 1}_{\{ \bold x (D) \leq n_1 \} }
 \sum_{i=1}^{m_1}
 \sum_{k=1}^{\bold x (D)}
 \beta^\mu(x_k) \cdot v_i (x_k)
 \partial_i
 f\left(\sum_{k=1}^{\bold x (D)}\varphi_1(x_k),\ldots,\sum_{k=1}^{\bold x (D)}\varphi_{m_1}(x_k)\right)
}
\\
 & &
 -
 {\bf 1}_{\{ \bold x (D) \leq n_1 \} }
 \sum_{i,j=1}^{m_1}
 \sum_{k=1}^{\bold x (D)} v_i (x_k)
 \sum_{l=1}^{\bold x (D)} v_j (x_l)
 \partial_i  \partial_j
 f\left(\sum_{k=1}^{\bold x (D)}\varphi_1(x_k),\ldots,\sum_{k=1}^{\bold x (D)}\varphi_{m_1}(x_k)\right)
\\
 & &
 +
 {\bf 1}_{\{ \bold x (D) \leq n_1 \} }
 \sum_{i=1}^{m_1}
 \sum_{k=1}^{\bold x (D)} \div v_i (x_k)
 \partial_i
 f\left(\sum_{k=1}^{\bold x (D)}\varphi_1(x_k),\ldots,\sum_{k=1}^{\bold x (D)}\varphi_{m_1}(x_k)\right)
\\
 & &
 +
 {\bf 1}_{\{ \bold x (D) \leq n_1 \} }
 \sum_{i=1}^{m_1}
 \sum_{k=1}^{\bold x (D)}
 U_{k,\bold x (D)} ( x_1,\ldots,x_{\bold x (D)} )
 \cdot
 v_i (x_k)
 \partial_i
 f\left(\sum_{k=1}^{\bold x (D)}\varphi_1(x_k),\ldots,\sum_{k=1}^{\bold x (D)}\varphi_{m_1}(x_k)\right)
\\
 & = &
 -\sum_{i=1}^{m_1}
 F_i (\bold{x} )
 \sum_{k=1}^{\bold x (D)}
 \beta^\mu(x_k)\cdot v_i (x_k)
\\
 & &
 -
 {\bf 1}_{\{ \bold x (D) \leq n_1 \} }
 \sum_{i,j=1}^{m_1}
 \sum_{k=1}^{\bold x (D)} v_i (x_k)
 \sum_{l=1}^{\bold x (D)} v_j (x_l)
 \partial_i  \partial_j
 f\left(\sum_{k=1}^{\bold x (D)}\varphi_1(x_k),\ldots,\sum_{k=1}^{\bold x (D)}\varphi_{m_1}(x_k)\right)
\\
 & &
 +
 \sum_{i=1}^{m_1}
 F_i (\bold{x} )
 \sum_{k=1}^{\bold x (D)} \div v_i (x_k)
\\
 & &
 +
 \sum_{i=1}^{m_1}
 F_i (\bold{x} )
 \nabla_{v_i}^{\mathrm{N}_{lf}} U[D] (\bold{x} )
,
\end{eqnarray*}
 which yields
\begin{align}
\mathcal{H}_DF(\bold{x})&=
 \sum_{i=1}^{m_1}\left(
 -\nabla_{v_i}^{\mathrm{N}_{lf}}F_i(\bold{x})
+(B_{v_i}^\mu(\bold{x})+\nabla_{v_i}^{\mathrm{N}_{lf}}U[D](\bold{x}))F_i(\bold{x})
 \right)
\nonumber\\
&=\sum_{i=1}^{m_1}\nabla_{v_i}^{\mathrm{N}_{lf}*} F_i(\bold{x}).\nonumber
\end{align}
 So, by Lemma~\ref{le:integrationbyparts} and
 since
 $\tilde{\cal S}_D \subset {\cal S}_D$,
 using obvious notation we have
\begin{align}
&\mathbb{E}[G(\bold{X}^D)\mathcal{H}_DF(\bold{X}^D)]=\sum_{i=1}^{m_1}\mathbb{E}\left[G(\bold{X}^D)\nabla_{v_i}^{\mathrm{N}_{lf}*}
F_i(\bold{X}^D)\right]=\sum_{i=1}^{m_1}\mathbb{E}\left[F_i(\bold{X}^D)\nabla_{v_i}^{\mathrm{N}_{lf}}G(\bold{X}^D)\right]\nonumber\\
&=\sum_{i=1}^{m_1}\mathbb{E}\left[F_i(\bold{X}^D)\sum_{l=1}^{\bold X (D)}\sum_{j=1}^{m_2}
\partial_j
g\left(\sum_{m=1}^{\bold X (D)}\gamma_1(X_m),\ldots,\sum_{m=1}^{\bold X (D)}\gamma_{m_2}(X_m)\right)\nabla\gamma_j(X_l)\cdot
\nabla\varphi_i(X_l)\right]\nonumber\\
&=\mathbb{E}\left[\sum_{l=1}^{\bold X (D)}\sum_{i=1}^{m_1}F_i(\bold{X}^D)\nabla\varphi_i(X_l)\cdot
\sum_{j=1}^{m_2}G_j(\bold{X}^D)\nabla\gamma_j(X_l)\right]\nonumber\\
&=\mathbb{E}\left[\sum_{i=1}^{\bold X (D)}\nabla_{X_i}^{\mathrm{N}_{lf}}F(\bold{X}^D)\cdot\nabla_{X_i}^{\mathrm{N}_{lf}}G(\bold{X}^D)\right].\nonumber
\end{align}
Finally, since $\mathcal{H}_D$ is symmetric and non-negative definite
the square root operator $\mathcal{H}_D^{1/2}$ is well-defined.
Relation~\eqref{eq:Dir2} follows by the properties of
$\mathcal{H}_D^{1/2}$.

\end{Proofy}

Let $\mathcal FC_b^\infty(D)$ denote the set of functionals $F:\mathrm N_f^D\to\R$ of the form
\begin{equation*}
F(\bold{x})=f\left(\sum_{k=1}^{\bold x (D)}\varphi_1(x_k),\ldots,\sum_{k=1}^{\bold x (D)}\varphi_m(x_k)\right),\qquad\bold x=\{x_1,\dots,x_{\mathbf x(D)}\}\in\mathrm N_f^D,
\end{equation*}
for some integer $m \geq 1$, $\varphi_1,\ldots,\varphi_m \in
{\cal C}^\infty(D)$, $f\in {\cal C}_b^\infty(\R^m)$.
For $F\in\mathcal FC_b^\infty(D)$, we naturally define
\begin{equation*}
\nabla_x^{\mathrm N_{lf}}F(\bold x):=\sum_{l=1}^m\partial_lf\left(\sum_{k=1}^{\bold x (D)}\varphi_1(x_k),\ldots,\sum_{k=1}^{\bold x (D)}\varphi_m(x_k)\right)\nabla\varphi_l(x),\qquad\bold x=\{x_1,\dots,x_{\mathbf x(D)}\}\in\mathrm N_f^D.
\end{equation*}
We conclude this section with the following proposition. 

\begin{proposition}
\label{prop:FCbinDomE}
Under $\bold{(H1)}-\bold{(H4)}$, we have that $\mathcal FC_b^\infty(D)\subset\mathrm{Dom}(\overline{\mathcal{H}_D^{1/2}})$ and
\begin{equation}
\label{eq:EDFG}
	\overline{\mathcal E}_D(F,G)=\mathbb{E}\left[\sum_{i=1}^{\bold X (D)}\nabla_{X_i}^{\mathrm{N}_{lf}}F(\bold{X}^D)\cdot\nabla_{X_i}^{\mathrm{N}_{lf}}G(\bold{X}^D)\right],\qquad F,G\in\mathcal FC_b^\infty(D).
\end{equation}
\end{proposition}
\begin{Proof}
For $F\in\mathcal FC_b^\infty(D)$, we clearly have $F^n(\bold x):=F(\bold x){\bf 1}_{\{|\bold x|(D)\le n\}}\in\tilde{\mathcal S}_D$, for any positive integer $n$.
By a straightforward computation, $F^m\to F$ in $\mathrm L^2_D$ as $n$ goes to infinity.
By the standard construction of the smallest closed extension of $\mathcal E_D$ (see e.g. \cite{fukushima}), to get that $F\in\mathrm{Dom}(\overline{\mathcal E}_D)$ and that $\overline{\mathcal E}_D(F,F)=\lim_{n\to\infty}\mathcal E_D(F^n,F^n)$, it suffices to prove that $\mathcal E_D(F^n-F^m,F^n-F^m)$ tends to zero as $m,n$ go to infinity.
This easily follows by the dominated convergence theorem.
Indeed for some positive $C>0$ and $m<n$,
\begin{align*}
\mathcal E_D(F^n-F^m,F^n-F^m)&=\mathbb E\left[\sum_{i=1}^{\bold X(D)}\|\nabla_{X_i}^{\mathrm N_{lf}}F(\bold X^D)\|^2{\bf 1}_{\{m<\bold X(D)\le n\}}\right]\\
&\le C\,\mathbb E\left[\bold X(D){\bf 1}_{\{m<\bold X(D)\le n\}}\right]
\xrightarrow[m,n\to\infty]{}0.
\end{align*}
Now, noticing that by the monotone convergence theorem,
\begin{equation*}
\mathcal E_D(F^n,F^n)=\mathbb E\left[\sum_{i=1}^{\bold X(D)}\|\nabla_{X_i}^{\mathrm N_{lf}}F(\bold X^D)\|^2{\bf 1}_{\{\bold X(D)\le n\}}\right]
\xrightarrow[n\to\infty]{}
\mathbb E\left[\sum_{i=1}^{\bold X(D)}\|\nabla_{X_i}^{\mathrm N_{lf}}F(\bold X^D)\|^2\right],
\end{equation*}
we have \eqref{eq:EDFG} with $G=F$ and we conclude by polarization.
\end{Proof}

\section{Diffusions associated to determinantal processes on $D\subset S$}\label{sec:diff}
 We start recalling some notions, see 
 Chapters IV and V in \cite{rockner}.
 Given $\pi$ in the set
 $\mathcal{P}(\mathrm{\ddot{N}}_f^D)$ of 
the probability measures on
$(\mathrm{\ddot{N}}_f^D,\mathcal{B}(\mathrm{\ddot{N}}_f^D))$,
 we call a $\pi$-stochastic process with state space $\mathrm{\ddot{N}}_f^D$
 the collection
$$
 \bold{M}_{D,\pi}=(\bold{\Omega},\mathcal{F},(\mathcal{F}_t)_{t\geq
0},(\bold{M}_t)_{t\geq
0},(\bold{\PP}_{\bold{x}})_{\bold{x}\in\mathrm{\ddot{N}}_f^D},\bold{\PP}_\pi)
$$
 where $\mathcal{F}:=\bigvee_{t\geq 0}\mathcal{F}_{t}$ is a $\sigma$-algebra
 on the set $\bold{\Omega}$,
 $(\mathcal{F}_t)_{t\geq 0}$ is the
${\bold \PP}_{\pi}$-completed filtration generated by
the process
$\bold{M}_t:\bold{\Omega} \longrightarrow \mathrm{\ddot{N}}_f^D$, 
$\bold{\PP}_{\bold{x}}$ is a
probability measure on $(\bold{\Omega},\mathcal{F})$
for all $\bold{x}\in\mathrm{\ddot{N}}_f^D$, and $\PP_\pi$ is
the probability measure on $(\bold{\Omega},\mathcal{F})$ defined
by
\[
\bold{\PP}_\pi( A ):=\int_{\mathrm{\ddot{N}}_f^D}\bold{\PP}_{\bold{x}}( A )\,\pi(\mathrm{d}\bold{x}),
\quad\text{$A \in\mathcal{F}$.}
\]
 A collection $(\bold{M}_{D,\pi},(\theta_t)_{t\geq 0})$ is called a
$\pi$-time homogeneous Markov process with state space
$\mathrm{\ddot{N}}_f^D$ if
$\theta_t:\bold{\Omega} \longrightarrow \bold{\Omega}$ is a shift operator, i.e.
$\bold{M}_s\circ\theta_t=\bold{M}_{s+t}$, $s,t \geq 0$;
 the map
$\bold{x}\mapsto \bold{\PP}_{\bold{x}}(A )$ is
measurable
for all $A \in\mathcal{F}$;
the time homogeneous Markov
property
$$\bold{\PP}_{\bold{x}}(\bold{M}_{t}\in
A \,|\,\mathcal{F}_s)=\bold{\PP}_{\bold{M}_s}(\bold{M}_{t-s}\in
A ),
 \qquad
 \bold{\PP}_{\bold{x}}-a.s.,
 \quad
 A \in\mathcal{F}, \quad
 0 \leq s \leq t,
 \quad
 \bold{x}\in\mathrm{\ddot{N}}_f^D
$$
 holds.
A $\pi$-time homogeneous Markov process
$(\bold{M}_{D,\pi},(\theta_t)_{t\geq 0})$ with state space
$\mathrm{\ddot{N}}_f^D$ is said to be $\pi$-tight on
$\mathrm{\ddot{N}}_f^D$ if $(\bold{M}_t)_{t\geq 0}$ is right-continuous with
left limits $\bold{\PP}_\pi$-almost surely;
$\bold{\PP}_{\bold{x}}(\bold{M}_0=\bold{x})=1$,
$\bold{x}\in\mathrm{\ddot{N}}_f^D$; the filtration
$(\mathcal{F}_t)_{t\geq 0}$ is right continuous; the following
strong Markov property holds:
$$\bold{\PP}_{\pi'}(\bold{M}_{t+\tau}\in
A
\,|\,\mathcal{F}_{\tau})=\bold{\PP}_{\bold{M}_\tau}(\bold{M}_{t}\in
A )$$ $\bold{\PP}_{\pi'}$-almost surely for all
$\mathcal{F}_t$-stopping time $\tau$,
$\pi'\in\mathcal{P}(\mathrm{\ddot{N}}_f^D)$, $A \in\mathcal{F}$
and $t\geq 0$, cf. Theorem IV.1.15 in \cite{rockner}. In addition,
a $\pi$-tight process on $\mathrm{\ddot{N}}_f^D$ is said to be a
$\pi$-special standard process on $\mathrm{\ddot{N}}_f^D$ if for
any $\pi'\in\mathcal{P}(\mathrm{\ddot{N}}_f^D)$ which is
equivalent to $\pi$ and all $\mathcal{F}_t$-stopping
times $\tau$, $(\tau_n)_{n\geq 1}$ such that $\tau_n\uparrow\tau$
then $\bold{M}_{\tau_n}$ converges to $\bold{M}_\tau$,
$\bold{\PP}_{\pi'}$-almost surely.
\\

\noindent
 In the following theorem, $\mathbb{E}_\bold{x}$ denotes the expectation under
 $\bold{\PP}_{\bold x}$, $\bold{x}\in\mathrm{\ddot{N}}_f^D$.
\begin{theorem}\label{thm:diff1}
 Assume $\bold{(H1)}-\bold{(H4)}$.
Then there exists a $\mathbb P_{D}$-tight special
standard process $(\bold{M}_{D,\mathbb P_{D}},(\theta_t)_{t\geq
0})$ on $\ddot{\mathrm N}_{f}^D$ such that:
\begin{enumerate}
\item
$\bold{M}_{D,\mathbb P_{D}}$ is a diffusion, in the sense that:
\begin{equation}\label{eq:diffusion}
\mathbb{P}_{\bold x}(\{\omega\in\bold\Omega:\,\,t\mapsto\bold{M}_t(\omega)\,\,\text{is
continuous on $[0,+\infty)$}\})=1,\qquad
\text{
$\overline{\mathcal{E}}_D$-a.e.
$\bold x\in\ddot{\mathrm N}_{f}^D$;} 
\end{equation}
\item
the transition semigroup of $\bold{M}_{D,\mathbb P_{D}}$ is given by
\begin{equation*}
p_t
F(\xi):=\mathbb{E}_{\bold x}[F(\bold{M}_t)],\quad\text{$\bold x\in\ddot{\mathrm N}_{f}^D$}, \quad
 F:\ddot{\mathrm N}_{f}^D \to \R
 \quad
 \mbox{square integrable},
\end{equation*}
and it
 is properly
associated with the Dirichlet form
$(\overline{\mathcal{E}}_D,\mathrm{Dom}(\overline{\mathcal{H}_D^{1/2}}))$
 in the sense that
 $p_t F$ is an $\overline{\mathcal{E}}_D$-a.c.,
 $\mathbb P_{D}$-version of
 $\exp(-t{\mathcal{H}_D^{gen}})F$,
 for all square integrable $F:\ddot{\mathrm N}_{f}^D \longrightarrow \R$ and $t>0$ (where ${\mathcal{H}_D^{gen}}$ is the generator of $\overline{\mathcal E}_D$);
\item
$\bold{M}_{D,\mathbb P_{D}}$ is unique up to
$\mathbb P_D$-equivalence (we refer the reader to Definition 6.3 page~140
in \cite{rockner} for the meaning of this notion);
\item
$\bold{M}_{D,\mathbb P_{D}}$ is $\mathbb P_D$-symmetric, i.e.
\begin{equation*}
	\mathbb E\left[{G}(\bold X^D)\,p_t{F}(\bold X^D)\right]=\mathbb E\left[{F}(\bold X^D)\,p_t{G}(\bold X^D)\right],
\end{equation*}
for $F,G\in\mathrm L^2_D$;
\item
$\bold{M}_{D,\mathbb P_{D}}$ has $\mathbb P_D$ as invariant measure.
\end{enumerate}
\end{theorem}
\begin{Proof}
We apply Theorem~4.13 of \cite[p. 308]{rockner1}.
Using the notation of \cite{rockner1}, set
\begin{equation*}
	\mathrm{S}(f,g)(x):=\nabla f(x)\cdot\nabla g(x),\qquad x\in\R^d,\ f,g\in\mathcal C^\infty(D),
\end{equation*}
where $\nabla$ is the usual gradient on $\R^d$, and
\begin{equation*}
	\mathrm{S}^\Gamma(F,G)(\bold x):=\sum_{i=1}^{\bold x(D)}\nabla_{x_i}^{\mathrm N_{lf}} F(\bold x)\cdot\nabla_{x_i}^{\mathrm N_{lf}} G(\bold x),\qquad\bold x=\{x_1,\dots,x_{\bold x(D)}\}\in\mathrm N_{f}^D,\ F,G\in
	\tilde{\mathcal S}_D.
\end{equation*}
Then, it is readily seen that $(\mathrm S,\mathcal C^\infty(D))$ satisfies conditions $(\mathcal D.1)$, $(\mathrm S.1)$, $(\mathrm S.2)$ and $(\mathrm S.3)$ of \cite{rockner1}, and $(\mathrm S^\Gamma,\tilde{\mathcal S}_D)$ satisfies condition $(\mathrm S^\Gamma.\mu)$ of \cite[p. 282]{rockner1}. Furthermore, $\mathbb P_D$ satisfies condition $(\mu.1)$ of \cite[p. 282]{rockner1} and condition $(Q)$ of \cite{rockner1} holds since $\R^d$ is complete (see Example~4.5.1 in \cite{rockner1}). The assumptions of Theorem~4.13 are therefore verified and the proof is completed.\end{Proof}
\subsubsection*{Non-collision property of the associated diffusions}
In the following, we will show the non-collision property of the diffusion constructed in the previous theorem which, roughly speaking, means that the diffusion takes values on $\mathrm N_{f}^D$.

We start by recalling the following lemma, which is borrowed from \cite{rocschm}.
\begin{lemma}
\label{lem:lemmarockner}
Assume $\bold{(H1)}-\bold{(H4)}$ and let $(\bold{M}_{D,\mathbb P_D},(\theta_t)_{t\geq 0})$ be the diffusion given by Theorem~\ref{thm:diff1}.
Let $u_n\in\mathrm{Dom}(\overline{\mathcal{E}}_D)$, $n\geq 1$, be such that: $u_n:\ddot{\mathrm{N}}_{f}^D\to\R$ is continuous, $u_n\to u$ point-wise in $\ddot{\mathrm N}_{f}^D$,
\begin{equation}
\label{eq:supDirichletfinite}
\sup_{n\geq 1}\overline{\mathcal{E}}_D(u_n,u_n)<\infty.
\end{equation}
Then $u$ is $\overline{\mathcal{E}}_D$-a.c. and, in particular,
$$
\bold{\PP}_{\bold{x}}(\{\omega\in\bold\Omega:\,\,t\mapsto u(\bold{M}_t)(\omega)\,\,\text{is
continuous on $[0,+\infty)$}\})=1,\qquad
\text{
$\overline{\mathcal E}_D$-a.e.
$\bold{x}\in\mathrm{N}_{f}^D$.}
$$
\end{lemma}
The next theorem provides the non-collision property. 
\begin{theorem}
\label{thm:noncollision}
 Assume $d \ge 2$, 
 and $\bold{(H1)}-\bold{(H4)}$.
Then 
\begin{equation*}
\bold{\PP}_{\bold{x}}(\{\omega\in\bold\Omega:\,\, \bold{M}_t(\omega)\in\mathrm{N}_{f}^D \quad\,\forall\, 0\le t < \infty\})=1,
\qquad\text{
$\overline{\mathcal E}_D$-a.e.
$\bold{x}\in\mathrm{N}_{f}^D$.}
\end{equation*}
\end{theorem}
\begin{Proof}
Since the proof is similar to the proof of Proposition $1$ in \cite{rocschm} we skip some details.
For every positive integer $a$, define
$u:=\bold{1}_N$, where
$$
N:=\{ \bold x \in \ddot{\mathrm{N}}_{f}^D\ : \ \sup_{x \in [-a,a]^d} {\bold x}(\{ x \}) \ge 2 \}.
$$
The claim follows if we prove that $u$ is $\overline{\mathcal{E}}_D$-a.c.. For this we are going to apply
Lemma \ref{lem:lemmarockner}. Define
$$
u_n (\bold x) = \Psi\left(\sup_{i \in A_n} \sum_{x \in \bold x} \phi_i(x) \right),\qquad\text{$n\geq 1$,}
$$
where $\Psi\in\mathcal C^\infty_b(\R)$ and $\phi_i\in\mathcal C^\infty(D)$ are chosen as in the proof of Proposition $1$ in \cite{rocschm}, and $A_n:= \mathbb{Z}^d \cap [-na,na]^d$.
Note that $u_n\in\mathrm{Dom}(\overline{\mathcal{H}_D^{1/2}})$ by Proposition~\ref{prop:FCbinDomE}. Furthermore, $u_n:\ddot{\mathrm{N}}_{f}^D\to\R$ is continuous and $u_n\to u$ point-wise by the proof of Proposition $1$ in \cite{rocschm}.
It remains to check \eqref{eq:supDirichletfinite}.
For $i=(i^{(1)},\dots,i^{(d)})\in \mathbb{Z}^d$ and $n\ge1$, we denote by $I_i^{(n)}$ the function defined by
$$
I_i^{(n)}(x):=\prod_{k = 1}^d {\bold 1}_{[-1/2, 3/2]}(n x^{(k)} - i^{(k)}),\qquad\text{$x = (x^{(1)},\dots,x^{(d)})\in D$.}
$$
As proved in  Proposition $1$ in \cite{rocschm}, the following upper bound holds:
\begin{equation}
\label{eq:boundonED}
\overline{\mathcal{E}}_D(u_n,u_n)
\le C n^2 \sum_{i \in A_n} \mathbb{E} \left[{\bold 1}_{\left\{ \sum_{j=1}^{\bold{X}(D)} I_i^{(n)}(X_j) \ge 2 \right\}}\sum_{j=1}^{\bold{X}(D)} I_i^{(n)}(X_j)\right],
\end{equation}
where $C>0$ is a positive constant.

Now, we upper-bound the r.h.s. of \eqref{eq:boundonED} by stochastic domination.
In \cite{georgiiyoo}, it is proved that
\begin{equation*}
	c[D](x,\bold x)\le J[D](x,x),\quad x\in D,\bold x\in\mathrm N_{f}^D,
\end{equation*} 
where we denote by $c[D](x,\bold x)$ the Papangelou conditional intensity of $\bold X^D$, see \cite{daley}.
Note that a Poisson process $\bold{Y}^D$ of mean measure $J[D](x,x)\mu(\mathrm dx)$ has Papangelou conditional intensity $J[D](x,x)$.
Therefore, by Theorem~1.1 in \cite{georgii}, we have that $\bold X^D$ is stochastically dominated by $\bold Y^D$, in the sense that
\begin{equation}\label{eq:stochasticdomination}
\mathbb{E}[f(\bold X^D)]\leq\mathbb{E}[f(\bold Y^D)]
\end{equation}
for any 
integrable $f:\mathrm{N}_{lf}\to\mathbb R$ such that
$f(\bold x) \le f(\bold y)$ whenever $\bold x \subset \bold y\in\mathrm N_{lf}$.
Since 
for any $i\in A_n$, the mapping $\bold x \mapsto {\bold 1}_{\left\{ \sum_{j=1}^{\bold{x}(D)} I_i^{(n)}(x_j) \ge 2 \right\}}\sum_{j=1}^{\bold{x}(D)} I_i^{(n)}(x_j)$ is
increasing, we then have 
\begin{equation}
\label{eq:bounde}
 \overline{\mathcal{E}}_D(u_n,u_n) \le C n^2 \sum_{i \in A_n} \mathbb{E}\left[{\bold 1}_{\left\{ \sum_{j=1}^{\bold{Y}(D)} I_j^{(n)}(Y_j) \ge 2 \right\}}\sum_{j=1}^{\bold{Y}(D)} I_i^{(n)}(Y_j)\right].
\end{equation}
By the properties of the Poisson process, the right hand side of \eqref{eq:bounde} is equal to
$$
C \, n^2 \sum_{i \in A_n}\left(1-\mathrm{e}^{-\int_DI_i^{(n)}(x)J[D](x,x)\,\mu(\mathrm{d}x)}\right)\int_DI_i^{(n)}(x)J[D](x,x)\,\mu(\mathrm{d}x),
$$
which is bounded above by
\begin{equation*}
C \, n^2 \sum_{i \in A_n}\left( \int_{D} I_i^{(n)}(x) J[D](x,x) \rho(x) \, \mathrm{d} x \right)^2.
\end{equation*}
By using the Cauchy-Schwarz inequality, this term is further bounded by
\begin{equation}\label{eq:ian}
C \, n^2 \sum_{i \in A_n}\left(\int_{D}I_i^{(n)}(x)\, \mathrm{d} x\right)\left(\int_{D}I_i^{(n)}(x) J[D](x,x)^2\rho(x)^2 \, \mathrm{d} x \right).
\end{equation}
We have on the one hand that
\[
\int_{\R^d}I_i^{(n)}(x)\,\mathrm{d}x=(2/n)^d
\]
and on the other hand for some constant $C'>0$,
\[
\int_{D}I_i^{(n)}(x) J[D](x,x)^2\rho(x)^2 \, \mathrm{d} x \le C'n^{-d},
\]
since $J$ and $\rho$ are bounded on $D$. 
Moreover, $\sharp(A_n)\leq\left(2 a n\right)^d$. Consequently, the quantity \eqref{eq:ian} is in turn bounded by
\begin{equation*}
C'' \, n^{2-d},\qquad\text{for some constant $C''>0$.}
\end{equation*}
The claim follows by the assumption
$d \ge 2$.
\end{Proof}
\section{An illustrating example}
\label{sec:examples}
 Let $S:=B(0,1)$ and $D:=B(0,R)\subset\R^2$ be the closed ball centered at the origin with radius
 $R \in (0,1)$ and $( \varphi_k^{(R)} )_{1\leq k\leq 3}$, the
 orthonormal subset of $\mathrm{L}^2(B(0,R),\ell )$ defined by
\[
\varphi_k^{(R)}(x):=
 \frac{1}{R}
 \sqrt{\frac{k+1}{\pi}}
 \left( \frac{x^{(1)}}{R}+i\frac{x^{(2)}}{R} \right)^k,
\quad x=(x^{(1)},x^{(2)})\in B(0,R),\ k=1,2,3,
\]
 where $\ell$ is the Lebesgue measure
 and $i:=\sqrt{-1}$ denotes the complex unit.
In this example, we consider the truncated Bergman kernel with at most $3$ points (see e.g. \cite{hough}) restricted to $D=B(0,R)$:
\[
K_{D}(x,y):=\sum_{k=1}^{3}R^{2(k+1)}\varphi_k^{(R)}(x)\overline{\varphi_k^{(R)}(y)},
 \qquad x,y \in D,
\]
and denote by $\mathcal{K}_{D}$ the associated integral operator,
 which is easily seen to be Hermitian and
trace class with non-zero eigenvalues
 $\kappa_k:=R^{2(k+1)}$, $k=1,2,3$.
 As a consequence, the spectrum of
 $\mathcal K_{D}$ is contained in $[0,1)$ and the triplet $(\mathcal
K_{D},K_{D},\ell)$ satisfies $(\bold{H1})$.
 In addition, $(\bold{H2})$ is trivially satisfied
 since the reference measure is the Lebesgue measure.
 The Janossy densities of $\bold{X}^D$ 
 defined in \eqref{def:janossy}
 are given by
\[
j_D^{n}(x_1,\ldots,x_n)=\mathrm{Det}(\bold{Id} - \mathcal{K}_{D})\,\mathrm{det}J[D](x_1,\dots,x_n), \qquad n = 1,2,3, \quad (x_1,\ldots,x_n)\in D^n,
\]
where the kernel $J[D]$ is given by
\[
J[D] (x,y):=\sum_{h=1}^{3}\frac{R^{2(h+1)}}{1-R^{2(h+1)}}\varphi_h^{(R)}(x)
\overline{\varphi_h^{(R)}(y)},
\]
 cf. (\ref{eq:jdecomp}). 
 Since $\bold X^D$ has at most $3$ points, see e.g. \cite{soshnikov}, we have $j_D^{n} = 0$, for $n \ge 4$.
To prove condition $(\bold{H3})$ it suffices to remark that the
function
\[
(x_1,\ldots,x_n)\to\mathrm{det}(J[D](x_p,x_q))_{1\leq p,q\leq
n}
\]
is continuously differentiable on $D^n$, for $n=1,2,3$.
To show that \text{\bf (H4)} is verified, we first consider the case $n=3$. Note that
\[
(J[D](x_p,x_q))_{1\leq p,q\leq
3}=A(x_1,x_2,x_3) A(x_1,x_2,x_3)^*,
\]
where the matrix $A(x_1,x_2,x_3):=(A_{ph} )_{1\leq p,h\leq 3}$ is given by
\[
A_{ph} :=\frac{R^{h+1}}{\sqrt{1-R^{2(h+1)}}}\varphi_h^{(R)}(x_p)
\]
and $A(x_1,x_2,x_3)^* $ denotes the transpose conjugate of $A(x_1,x_2,x_3) $.
Hence,
\[
\mathrm{det}\,J[D](x_1,x_2,x_3)=|\mathrm{det}\, A(x_1,x_2,x_3) |^2,
\]
and since the previous determinant can be rewritten involving a Vandermonde determinant, we have
\[
\mathrm{det}\, A(x_1,x_2,x_3)= \prod_{p=1}^3 \sqrt{\frac{1+p}{\pi(1-R^{2(p+1)})}}
\left(
 \prod_{p=1}^{3}(x_p^{(1)}+ix_p^{(2)})
\right)
 \prod_{1\leq p<q \leq 3}
 ((x_p^{(1)}-x_q^{(1)})+i(x_p^{(2)}-x_q^{(2)})).
\]
\\

\noindent
 Note that \text{\bf (H4)} with $n=3$ is exactly
\begin{equation*}
\int_{
 D^3}
 \left|
 \frac{\partial_{x_i^{(h)}}
|\mathrm{det}\, A(x_1,x_2,x_3)|^2
 \partial_{x_j^{(k)}}
| \mathrm{det}\,A(x_1,x_2,x_3)|^2
 }{
 |\mathrm{det}\,A(x_1,x_2,x_3)|^2
 }
 \right|
 {\bf1}_{ \{|\mathrm{det}\,A(x_1,x_2,x_3)|>0\}}
\,\mathrm{d}x_1\mathrm{d}x_2\mathrm{d}x_3<\infty,
\end{equation*}
 for all $1\leq i,j\leq 3$ and $1\leq h,k\leq 2$,
 and so it suffices to check
$$
\int_{
 B(0,R)^3}
 \left|
 \frac{\partial_{x_1^{(1)}}
|\mathrm{det}\, A(x_1,x_2,x_3)|^2
 }{
 |\mathrm{det}\,A(x_1,x_2,x_3)|^2
 }
 \right|
  {\bf1}_{ \{|\mathrm{det}\,A(x_1,x_2,x_3)|>0\}}
\,\mathrm{d}x_1\mathrm{d}x_2\mathrm{d}x_3<\infty.
$$
This latter integral reduces to
$$
\int_{
 B(0,R)^3}
 \left|
 \frac{2 x_1^{(1)}
 }{
(x_1^{(1)})^2 + (x_1^{(2)})^2
 }
 +
 2 \sum_{j=2}^3
  \frac{x_1^{(1)} - x_j^{(1)}
 }{
(x_1^{(1)} - x_j^{(1)})^2 + (x_1^{(2)}- x_j^{(2)})^2
 }
 \right|
\,\mathrm{d}x_1\mathrm{d}x_3\mathrm{d}x_3,
$$
which is indeed finite.
Consequently, we proved that \text{\bf (H4)} is verified for $n \ge 3$ (indeed it is trivially satisfied for $n > 3$).
Now, consider $n=1,2$. We have again
\[
J[D](x_1,\dots,x_n)=A(x_1,\dots,x_n)A(x_1,\dots,x_n)^*,
\]
where this time, $A(x_1,\dots,x_n) $ is a rectangular $n\times3$ matrix given by
 $A(x_1,\ldots,x_n):=(A_{ph} )_{1\leq p\le n,1\le h\leq 3}$ is given by
\[
A_{ph} :=\frac{R^{h+1}}{\sqrt{1-R^{2(h+1)}}}\varphi_h^{(R)}(x_p).
\]
Recall the Cauchy-Binet formula:
\begin{equation}
\label{eq:eq1examples}
\mathrm{det}\,J[D](x_1,\dots,x_n)=\sum_{1 \le i_1 < i_2 < \dots < i_n \le 3} | \mathrm{det}\, A^{i_1,\dots,i_n}(x_1,\dots,x_n)|^2,
\end{equation}
where 
\[
A^{i_1,\dots,i_n}_{ph} :=\frac{R^{i_h+1}}{\sqrt{1-R^{2(i_h+1)}}}\varphi_{i_h}^{(R)}(x_p),
\qquad 1\le p,h \le n,
\]
defines a square matrix. We now consider fixed $1 \le i_1 < i_2 < \dots < i_n \le 3$ and evaluate $| \mathrm{det}\, A^{i_1,\dots,i_n}(x_1,\dots,x_n)|^2$. We note that
\begin{equation}
\label{eq:eq2examples}
| \mathrm{det}\, A^{i_1,\dots,i_n}(x_1,\dots,x_n)|^2 = \prod_{p=1}^n \frac{1+i_p}{\pi(1-R^{2(i_p+1)})}
\left|
V_{i_1,\dots,i_n}(x_1,\dots,x_n)
\right|^2,
\end{equation}
where
\[
V_{i_1,\dots,i_n}(x_1,\dots,x_n) := \mathrm{det}\left( \left(x_h^{i_p}\right)_{1\le p,h \le n} \right)
\]
is known in the literature as the generalized Vandermonde determinant.
By definition, the Schur polynomial $s_\lambda$ is the ratio between the generalized Vandermonde determinant and the classical Vandermonde determinant. 
More precisely,
\begin{equation}
\label{eq:eq3examples}
\mathrm{det}\left( \left(x_h^{i_p}\right)_{1\le p,h \le n} \right) = \mathrm{det}\left( \left(x_h^{p-1}\right)_{1\le p,h \le n} \right)s_{\lambda(i_1,\dots,i_n)}(x_1,\dots,x_n),
\end{equation}
where ${\lambda(i_1,\dots,i_n)}:= (i_n-n+1,\dots,i_2-1,i_1)$.
Combining \eqref{eq:eq1examples}, \eqref{eq:eq2examples} and \eqref{eq:eq3examples}, we have
\begin{multline}
\label{eq:detjbergman}
\mathrm{det}\,J[D](x_1,\dots,x_n) = \left| \prod_{1\leq p<q \leq n}
 ((x_p^{(1)}-x_q^{(1)})+i(x_p^{(2)}-x_q^{(2)}))\right|^2
\\
\nonumber
 \sum_{1 \le i_1 < i_2 < \dots < i_n \le 3} \left( \prod_{p=1}^n \frac{1+i_p}{\pi(1-R^{2(i_p+1)})}\right)
 | s_{\lambda(i_1,\dots,i_n)}(x_1,\dots,x_n)|^2.
\end{multline}
\noindent
For $n=1$, one has $s_{\lambda(1)}(x)=x$, $s_{\lambda(2)}(x)=x^2$ and $s_{\lambda(3)}(x)=x^3$, see e.g. \cite{heineman}, and therefore
\begin{equation*}
\mathrm{det}\,J[D](x) = \frac{2}{\pi(1-R^{4})}
 |x |^2+ \frac{3}{\pi(1-R^{6})} |x|^4+\frac{4}{\pi(1-R^{8})} |x|^6.
\end{equation*}
Thus,
\begin{align*}
	&\partial_{x^{(p)}}\ln\left(\mathrm{det}\,J[D](x)\right)
	=\frac{2}{\pi(1-R^{4})}\frac{2x^{(p)}}{\frac{2}{\pi(1-R^{4})}
 |x |^2+ \frac{3}{\pi(1-R^{6})} |x|^4+\frac{4}{\pi(1-R^{8})} |x|^6}\\
 &\qquad\qquad\qquad\qquad\qquad\qquad\qquad\qquad+\frac{3}{\pi(1-R^{6})}\frac{4x^{(p)}|x|^2}{\frac{2}{\pi(1-R^{4})}|x |^2+ \frac{3}{\pi(1-R^{6})} |x|^4+\frac{4}{\pi(1-R^{8})} |x|^6}
 \\
 &\qquad\qquad\qquad\qquad\qquad\qquad\qquad\qquad+\frac{4}{\pi(1-R^{8})}\frac{6x^{(p)}|x|^4}{\frac{2}{\pi(1-R^{4})}
 |x |^2+ \frac{3}{\pi(1-R^{6})} |x|^4+\frac{4}{\pi(1-R^{8})} |x|^6},
\end{align*}
and therefore
\begin{equation*}
	\left|\partial_{x^{(p)}}\ln\left(\mathrm{det}\,J[D](x)\right)\right|
	\le\frac{2\left|x^{(p)}\right|}{ |x |^2}
 +\frac{4\left|x^{(p)}\right|}{ |x|^2}
 +\frac{6\left|x^{(p)}\right|}{ |x|^2},
\end{equation*}
which is integrable on $D$.
For $n=2$, one has $s_{\lambda(1,2)}(x,y)=1$, $s_{\lambda(1,3)}(x,y)=x+y$ and $s_{\lambda(2,3)}(x,y)=xy$, see e.g. \cite{heineman}, and therefore
\begin{multline*}
\mathrm{det}\,J[D](x,y) = \left|x-y\right|^2\\
\times\left(\frac{2}{\pi(1-R^{4})}\frac{3}{\pi(1-R^{6})}+ \frac{2}{\pi(1-R^{4})}\frac{4}{\pi(1-R^{8})} |x+y|^2+\frac{3}{\pi(1-R^{6})}\frac{4}{\pi(1-R^{8})} |xy|^2\right).
\end{multline*}
Note that the differential of the logarithm of $\left|x-y\right|^2$ gives rise to a locally integrable term.
So it remains to check that the differential of the logarithm of the second term, hereafter denoted by $\Psi_D(x,y)$, is integrable on $D^2$.
By symmetry of the Schur polynomials, it suffices to check that the derivative of $\Psi_D(x,y)$ with respect to $x^{(p)}$ is integrable on $D^2$.
We have
\begin{equation*}
	\left|\partial_{x^{(p)}}\Psi_D(x,y)\right|
	\le\frac{2\left|x^{(p)}+y^{(p)}\right|}{|x+y|^2}+\frac{2\left|x^{(p)}\right||y|^2}{|xy|^2}
	=\frac{2\left|x^{(p)}+y^{(p)}\right|}{|x+y|^2}+\frac{2x^{(p)}}{|x|^2},
\end{equation*}
and the claim follows by noticing that the r.h.s. is integrable on $D^2$.

\footnotesize

\def\cprime{$'$} \def\polhk#1{\setbox0=\hbox{#1}{\ooalign{\hidewidth
  \lower1.5ex\hbox{`}\hidewidth\crcr\unhbox0}}}
  \def\polhk#1{\setbox0=\hbox{#1}{\ooalign{\hidewidth
  \lower1.5ex\hbox{`}\hidewidth\crcr\unhbox0}}} \def\cprime{$'$}

\end{document}